\documentclass[final]{siamltex}
\usepackage{amsmath,amsfonts,amssymb}
\newtheorem{teor}{Theorem}[section]
\newtheorem{lema}[teor]{Lemma}
\newtheorem{prop}[teor]{Proposition}

\newtheorem{defi}[teor]{Definition}
\newtheorem{nota}[teor]{Remark}
\numberwithin{equation}{section}
\newcommand{\R}{\mathbb{R}}

\newcommand{\N}{\mathbb{N}}

\newcommand{\Om}{\Omega}
\newcommand{\w}{\omega}

\newcommand{\ep}{\varepsilon}

\newcommand{\di}{\textsf{d}}
\newcommand{\wt}{\widetilde}
 \DeclareMathOperator{\car}{card}
\newcommand{\n}[1]{\| #1 \|}

\title{Neutral functional differential equations with applications to
        compartmental systems\thanks{The authors were partly supported
        by Junta de Castilla y Le\'{o}n under project VA024A06, and
        M.E.C. under project MTM2005-02144.}}

\author{V{\'\i}ctor Mu\~noz-Villarragut, Sylvia Novo
        \and Rafael Obaya\thanks{Departamento de Matem\'{a}tica
        Aplicada, E.T.S. de Ingenieros Industriales, Universidad de
        Valladolid, 47011 Valladolid, Spain ({\tt vicmun@wmatem.eis.uva.es}},
        {\tt sylnov@wmatem.eis.uva.es},
        {\tt rafoba@wmatem.eis.uva.es}).}

\begin{document}

\vspace*{-1.5cm}

\noindent\begin{tabular}{|l|}
  \hline
  V. Muñoz-Villarragut, S. Novo, R. Obaya,\\
  Neutral Functional Differential Equations with
  Applications to Compartmental Systems,\\
  SIAM Journal on Mathematical Analysis, Volume 40,
  Number 3, 2008, Pages 1003-1028.\\
  https://doi.org/10.1137/070711177\\
  \copyright ~Society for Industrial and Applied Mathematics\\
  \hline
\end{tabular}

\vspace{1cm}

\maketitle

\begin{abstract}
We study the monotone skew-product semiflow generated by a family
of neutral functional differential equations with infinite delay
and stable $D$-operator. The stability properties of $D$ allow us
to introduce a new order and to take the neutral family to a
family of functional differential equations with infinite delay.
Next, we establish the 1-covering property of omega-limit sets
under the componentwise separating property and uniform stability.
Finally, the obtained results are applied to the study of the
long-term behavior of the amount of material within the
compartments of a neutral compartmental system with infinite
delay.
\end{abstract}

\begin{keywords}
Non-autonomous dynamical systems, monotone skew-product semiflows,
neutral functional differential equations, infinite delay,
compartmental systems
\end{keywords}

\begin{AMS}
37B55, 34K40, 34K14
\end{AMS}

\pagestyle{myheadings} \thispagestyle{plain} \markboth{V.
MU\~NOZ-VILLARRAGUT, S. NOVO AND R. OBAYA}{NEUTRAL FUNCTIONAL
DIFFERENTIAL EQUATIONS}

\section{Introduction}\label{intro}
After the pioneering work of Hale and Meyer~\cite{hame}, the
theory of neutral functional differential equations (NFDEs for
short) aroused a considerable interest and a fast development
ensued. At present a wide collection of theoretical and practical
results make up the main body of the theory of NFDEs (see
Hale~\cite{hale}, Hale and Verduyn Lunel~\cite{hale2},
Kolmanovskii and  Myshkis~\cite{kolm} and Salamon~\cite{sala},
among many others). In particular, a substantial number of results
for delayed functional differential equations (FDEs for short)
have been generalized for NFDEs solving new and challenging
problems in these extensions.
\par
In this paper we provide a dynamical theory for non-autonomous
monotone NFDEs with infinite delay and autonomous stable
$D$-operator in the line of the results by Jiang and
Zhao~\cite{jizh} and Novo {\it et al.}~\cite{NOS2007}. We assume
some recurrence properties on the temporal variation of the NFDE.
Thus, its solutions induce  a skew-product semiflow with minimal
flow on the base. In particular, the uniform almost periodic and
almost automorphic cases are included in this formulation. The
skew-product formalism permits the analysis of the dynamical
properties of the trajectories using methods of ergodic theory and
topological dynamics.\par Novo {\it et al.}~\cite{NOS2007} study
the existence of recurrent solutions of non-autonomous  FDEs with
infinite delay using the phase space $BU\subset
C((-\infty,0],\R^m)$ of bounded and uniformly continuous functions
with the supremum norm. Assuming some technical conditions on the
vector field, it is shown that every bounded solution is
relatively compact for the compact-open topology and its
omega-limit set admits a flow extension. An alternative method for
the study of recurrent solutions of almost periodic FDEs with
infinite delay makes use of a {\em fading memory\/} Banach phase
space (see Hino {\it et al.}~\cite{hino} for an axiomatic
definition and main properties). Since this kind of spaces contain
$BU$ and, under natural assumptions, the restriction of the norm
topology to the closure of a bounded solution agrees with the
compact-open topology, it seems that the approach considered in
Novo {\it et al.}~\cite{NOS2007} becomes natural in many cases of
interest.
\par
In this paper we consider NFDEs with linear
autonomous operator $D$ defined on $BU$ which is continuous for
the norm, continuous for the compact-open topology on bounded sets
and atomic at zero. We obtain an integral representation of $D$ by
means of Riesz theorem $Dx=\int_{-\infty}^0 [\,d\mu] \,x$, where
$\mu$ is a real Borel measure with finite total variation. The
convolution operator $\widehat D$ defined by $\widehat
Dx(s)=\int_{-\infty}^0 [d\mu(\theta)] \,x(\theta+s)$ maps $BU$
into $BU$. We prove that if $D$ is stable in the sense of
Hale~\cite{hale}, then $\widehat D$ is an isomorphism of $BU$
which is continuous for the norm and continuous for the
compact-open topology on bounded sets.
 Moreover, $\widehat D^{-1}$ inherits these same properties and it is
associated to a linear stable operator $D^*$. In fact, the
mentioned behavior of $\widehat D^{-1}$ characterizes the
stability of the operator $D$. The proofs are self-contained and
only require quantitative estimates associated to the stability of
the operator $D$.
\par
Staffans~\cite{staf83} shows that every NFDE
with finite delay and autonomous stable $D$-operator can be
written as  a FDE with infinite delay in an appropriate fading
memory space. A more systematic study on the inversion of the
convolution operator $\widehat D$ can be found in Gripenberg  {\it
et al.}~\cite{GLS1990}.
 The papers by Haddock {\it et
al.}~\cite{HKW1990, HKW1994}, Arino and Bourad~\cite{AB1990},
among others, make a systematic use of these ideas which have a
theoretical and practical interest. We give a version of the above
results for infinite delay NFDEs. It is obvious that the inversion
of the convolution operator $\widehat D$ on $BU$ allows us to
transform the original equation into a retarded non-autonomous FDE
with infinite delay. In addition, we transfer the dynamical theory
of Jiang and Zhao~\cite{jizh} and Novo {\it et al.}~\cite{NOS2007}
to non-autonomous monotone NFDEs with infinite delay and
autonomous stable $D$-operator. In an appropriate dynamical
framework we assume that the trajectories are bounded, uniformly
stable and satisfy a componentwise separating property, and we
show that the omega-limit sets are all copies of the base. It is
important to mention that no conditions of strong monotonicity are
required, which permits the application of the results under
natural physical conditions.
\par
In this paper we provide a detailed description of the long-term
behavior of the dynamics in some classes of compartmental systems
extensively studied in the literature. Compartmental systems have
been used as mathematical models for the study of the dynamical
behavior of many processes in biological and physical sciences,
which depend on local mass balance conditions (see
Jacquez~\cite{JJ1996}, Jacquez and Simon~\cite{JS1993,JS2002}, and
references therein). Some initial results for models described by
FDEs with finite and infinite delay can be found in
Gy\"{o}ri~\cite{gyori} and Gy\"{o}ri and Eller~\cite{gyell}. The
paper by Arino and Haourigui~\cite{ARI1987} proves the existence
of almost periodic solutions for compartmental systems described
by almost periodic finite delay FDEs. NFDEs represent systems
where the compartments produce or swallow material. Gy\"{o}ri and
Wu~\cite{gyoriwu} and Wu and Freedman~\cite{WF1991} study
autonomous NFDEs with finite and infinite delay similar to those
considered in this paper. We provide a non-autonomous version,
under more general assumptions, of the monotone theory for NFDEs
included in Wu and Freedman~\cite{WF1991} and Wu~\cite{jwu1991}.
More precise results for the case of scalar NFDEs can be found in
Arino and Bourad~\cite{AB1990} and Krisztin and Wu~\cite{KW1996}.
\par
We study the dynamics of monotone compartmental systems  in terms
of the geometrical structure of the pipes connecting the
compartments. Irreducible subsets of the
 set of indices detect the occurrence of subsystems on the
complete system which reduce the dimension of the problem to
study. When the system is closed, the total mass is an invariant
continuous function which implies the stability and boundedness of
the solutions. In particular,  the omega-limit set of every
solution is a minimal set and a copy of the base. In a general
compartmental system the existence of a bounded solution assures
that every solution is bounded and uniformly stable. We first
check that when there is no inflow of material then all the
compartments of each irreducible subset with some outflow of
material are empty on minimal subsets. On the contrary, when the
solutions remain bounded and there is inflow of material in some
compartment of an irreducible subset, then all the minimal sets
agree  for the indices of this irreducible subset and all the
compartments contain some material. Finally, we describe natural
physical conditions which ensure the existence of a unique minimal
set asymptotically stable.
\par
This paper is organized as follows. Basic notions in topological
dynamics, used along the rest of the sections, are stated in
Section~\ref{prelim}. Section~\ref{stableD} is devoted to the
study of general and stability properties of linear autonomous
operators from $BU$ to $\R^m$, as well as the behavior of
solutions of the corresponding homogeneous and nonhomogeneous
equations given by them. In Section~\ref{monotone}, we study the
monotone skew-product semiflow generated by a family of NFDEs with
infinite delay and stable $D$-operator. In particular, we
establish the 1-covering property of omega-limit sets under the
componentwise separating property and uniform stability. These
results are applied in Section~\ref{compartmental} to show that
the solutions of a compartmental system given by a monotone NFDE
with infinite delay are asymptotically of the same type as the
transport functions. Finally, Section~\ref{long-term} deals with
the long-term behavior of its solutions in terms of the
geometrical structure of the pipes, as explained above.
\section{Some preliminaries}\label{prelim}
Let $(\Om,d)$ be a compact metric space. A real {\em continuous
flow \/}  $(\Om,\sigma,\R)$ is defined by a continuous mapping
$\sigma: \R\times \Om \to  \Om,\; (t,\w)\mapsto \sigma(t,\w)$
satisfying
\begin{enumerate}
\renewcommand{\labelenumi}{(\roman{enumi})}
\item $\sigma_0=\text{Id},$
\item $\sigma_{t+s}=\sigma_t\circ\sigma_s$ for each $s$, $t\in\R$,
\end{enumerate}
where $\sigma_t(\w)=\sigma(t,\w)$ for all $\w \in \Om$ and $t\in
\R$. The set $\{ \sigma_t(\w) \mid t\in\R\}$ is called the {\em
orbit\/} or the {\em trajectory\/} of the point $\w$. We say that
a subset $\Om_1\subset \Om$ is {\em $\sigma$-invariant\/} if
$\sigma_t(\Om_1)=\Om_1$ for every $t\in\R$.  A subset
$\Om_1\subset  \Om$  is called {\em minimal \/} if it is compact,
$\sigma$-invariant and its only nonempty compact
$\sigma$-invariant subset is itself. Every compact and
$\sigma$-invariant  set contains a minimal subset; in particular
it is  easy to prove that a compact $\sigma$-invariant subset is
minimal if and only if every trajectory is dense. We say that the
continuous flow $(\Om,\sigma,\R)$ is {\em recurrent\/} or {\em
minimal\/} if $\Om$ is minimal.
\par
The flow $(\Om,\sigma,\R)$ is {\em distal\/} if for any two
distinct points $\w_1,\,\w_2\in\Om$ the orbits keep at a positive
distance, that is, $\inf_{t\in
\R}d(\sigma(t,\w_1),\sigma(t,\w_2))>0$. The flow $(\Om,\sigma,\R)$
is {\em almost periodic\/} when for every $\varepsilon
> 0 $ there is a $\delta >0$ such that, if $\w_1$, $\w_2\in\Om$
with $d(\w_1,\w_2)<\delta$, then
$d(\sigma(t,\w_1),\sigma(t,\w_2))<\varepsilon$ for every $t\in
\R$. If $(\Om,\sigma,\R)$ is almost periodic, it is distal. The
converse is not true; even if $(\Om,\sigma,\R)$ is minimal and
distal, it does not need to be almost periodic. For the basic
properties of almost periodic and distal flows we refer the reader
to Ellis~\cite{elli} and Sacker and Sell~\cite{sase1}.
\par
A {\em flow homomorphism\/} from another continuous flow
$(Y,\Psi,\R)$ to $(\Om,\sigma,\R)$ is a continuous map $\pi\colon
Y\to \Om$ such that $\pi(\Psi(t,y))=\sigma(t,\pi(y))$ for every
$y\in Y$ and $t\in\R$. If $\pi$ is also bijective, it is called a
{\em flow isomorphism\/}. Let $\pi:Y \to \Om$ be a surjective flow
homomorphism and suppose $(Y,\Psi,\R)$ is minimal (then, so is
$(\Om,\sigma,\R)$). $(Y,\Psi,\R)$ is said to be an {\em almost
automorphic extension\/} of $(\Om,\sigma,\R)$ if there is $\w\in
\Om$ such that $\car (\pi^{-1}(\w))=1$. Then, actually $\car
(\pi^{-1}(\w))=1$ for $\w$ in a residual subset $\Om_0\subseteq
\Om$; in the nontrivial case $\Om_0\subsetneq \Om$ the dynamics
can be very complicated. A minimal flow $(Y,\Psi,\R)$ is {\em
almost automorphic\/} if it is an almost automorphic extension of
an almost periodic minimal flow $(\Om,\sigma,\R)$. We refer the
reader to the work of Shen and Yi~\cite{shyi} for a survey of
almost periodic and almost automorphic dynamics.
\par
Let $E$ be a complete metric space and  $\R^+=\{t\in\R\,|\,t\geq
0\}$. A {\em semiflow} $(E,\Phi,\R^+)$ is determined by a
continuous map $\Phi: \R^+\times E \to E,\; (t,x)\mapsto
\Phi(t,x)$ which satisfies
\begin{enumerate}
\renewcommand{\labelenumi}{(\roman{enumi})}
\item $\Phi_0=\text{Id},$
\item $\Phi_{t+s}=\Phi_t \circ \Phi_s\;$ for all  $\; t$, $s\in\R^+,$
\end{enumerate}
where $\Phi_t(x)=\Phi(t,x)$ for each $x \in E$ and $t\in \R^+$.
The set $\{ \Phi_t(x)\mid t\geq 0\}$ is the {\em semiorbit\/} of
the point $x$. A subset  $E_1$ of $E$ is {\em positively
invariant\/} (or just $\Phi$-{\em invariant\/}) if
$\Phi_t(E_1)\subset E_1$ for all $t\geq 0$. A semiflow
$(E,\Phi,\R^+)$ admits a {\em flow extension\/} if there exists a
continuous flow  $(E,\wt \Phi,\R)$ such that $\wt
\Phi(t,x)=\Phi(t,x)$ for all $x\in E$ and $t\in\R^+$. A compact
and positively invariant subset admits a flow extension if the
semiflow restricted to it admits one.
\par
Write $\R^-=\{t\in\R\,|\,t\leq 0\}$. A {\em backward orbit\/} of a
point  $x\in E$ in the semiflow $(E,\Phi,\R^+)$ is a continuous
map  $\psi:\R^-\to E$ such that $\psi(0)=x$ and for each $s\leq 0$
it holds that  $\Phi(t,\psi(s))=\psi(s+t)$ whenever $0\leq t\leq
-s$. If for $x\in E$ the semiorbit $\{\Phi(t,x)\mid t\ge 0\}$ is
relatively compact, we can consider the {\em omega-limit set\/} of
$x$,
\[
\mathcal{O}(x)=\bigcap_{s\ge 0}{\rm closure}{\{\Phi(t+s,x)\mid
t\ge 0\}}\,,
\]
which is a nonempty compact connected and $\Phi$-invariant set.
Namely, it consists of the points $y\in E$ such that $y=\lim_{n\to
\infty} \Phi(t_n,x)$ for some sequence $t_n\uparrow \infty$. It is
well-known that every $y\in\mathcal{O}(x)$ admits a backward orbit
inside this set. Actually, a compact positively invariant set $M$
admits a flow extension if every point in $M$ admits a unique
backward orbit which remains inside the set $M$ (see Shen and
Yi~\cite{shyi},~part~II).
\par
A compact positively invariant set $M$ for the semiflow
$(E,\Phi,\R^+)$ is {\em minimal\/} if it does not contain any
other nonempty compact positively invariant set than itself. If
$E$ is minimal, we say that the semiflow is minimal.
\par
A semiflow is {\em of skew-product type\/} when it is defined on a
vector bundle and has a triangular structure; more precisely, a
semiflow $(\Om\times X,\tau,\,\R^+)$ is a {\em skew-product\/}
semiflow over the product space $\Om\times X$, for a compact
metric space $(\Om,d)$ and a complete metric space
$(X,\textsf{d})$, if the continuous map $\tau$ is as follows:
\begin{equation}\label{skewp}
 \begin{array}{cccl}
 \tau \colon  &\R^+\times\Om\times X& \longrightarrow & \Om\times X \\
& (t,\w,x) & \mapsto &(\w{\cdot}t,u(t,\w,x))\,,
\end{array}
\end{equation}
where $(\Om,\sigma,\R)$ is a  real continuous flow
$\sigma:\R\times\Om\rightarrow\Om$, $\,(t,\w)\mapsto \w{\cdot}t$, called
the  {\em base flow\/}. The skew-product semiflow~\eqref{skewp} is
{\em linear\/} if $u(t,\w,x)$ is linear in $x$ for each
$(t,\w)\in\R^+\times\Om$.
\par
Now, we introduce some definitions concerning the stability of the
trajectories. A forward orbit $\{\tau(t,\w_0,x_0)\,|\;t\geq 0\}$
of the skew-product semiflow~\eqref{skewp} is said to be {\em
uniformly stable\/} if for every $\varepsilon>0$ there is a
$\delta(\varepsilon)>0$, called the {\em modulus of uniform
stability\/}, such that, if $s\geq 0$ and
$\textsf{d}(u(s,\w_0,x_0),x)\leq \delta(\ep)$ for certain $x\in
X$, then for each $t\geq 0$,
\[
\di(u(t+s,\w_0,x_0),u(t,\w_0{\cdot}s,x))=\di(u(t,\w_0{\cdot}s,u(s,\w_0,x_0)),u(t,\w_0{\cdot}s,x))\leq
\varepsilon \,.
\]
A forward orbit $\{\tau(t,\w_0,x_0)\,|\;t\geq 0\}$ of the
skew-product semiflow~\eqref{skewp} is said to be {\em uniformly
asymptotically stable\/} if it is uniformly stable and there is a
$\delta_0>0$ with the following property: for each $\varepsilon>0$
there is a $t_0(\varepsilon)>0$ such that, if $s\geq 0$ and
$\di(u(s,\w_0,x_0),x)\leq \delta_0$, then
\[\di(u(t+s,\w_0,x_0),u(t,\w_0{\cdot}s,x))\leq \varepsilon \quad \text{ for
each } t\geq t_0( \varepsilon)\,.\]
\section{Stable $D$-operators}\label{stableD}
 We consider the Fr\'{e}chet space
$X=C((-\infty,0],\R^m)$ endowed with the compact-open topology,
i.e.~the topology of uniform convergence over compact subsets,
which is a metric space for the distance
\[\textsf{d}(x,y)=\sum_{n=1}^\infty \frac{1}{2^n}\frac
{\n{x-y}_n}{1+\n{x-y}_n}\,,\quad x,y\in X\,,
\]
where $\n{x}_n=\sup_{s\in[-n,0]}\n{x(s)}$, and $\n{\cdot}$ denotes the
maximum norm on $\R^m$.
\par
Let $BU\subset X$ be the Banach space
\[BU=\{x\in X\mid  x \text{ is bounded and
  uniformly continuous}\}\]
 with the supremum norm  $\n{x}_\infty=\sup_{s\in(-\infty,0]} \n{x(s)}$.
Given $r>0$,  we will denote
\[B_r=\{x\in BU \mid \n{x}_\infty \leq r\}\,.\] As usual, given
$I=(-\infty,a]\subset\R$,  $t\in I$ and a continuous function
$x:I\to\R^m$, $x_t$ will denote the element of $X$ defined by
$x_t(s)=x(t+s)$ for $s\in (-\infty,0]$.
\par
This section is devoted to the study of general and stability
properties of linear autonomous operators $D\colon BU\to \R^m$, as
well as the behavior of solutions of the corresponding homogeneous
equation $Dx_t=0$, $t\geq 0$, and nonhomogeneous equations
$Dx_t=h(t)$, $t\geq 0$ for $h\in C( [0,\infty),\R^m)$. We will
assume that:\smallskip
\begin{itemize}
\item[(D1)] $D$ is linear and continuous for the norm.
\item[(D2)] For each $r>0$, $D\colon B_r\to \R^m$ is continuous when we take
the restriction of the compact-open topology to $B_r$, i.e. if
$x_n\stackrel{\textsf{d}\;}\to x$ as $n\to\infty$ with $x_n$,
$x\in B_r$, then $\lim_{\,n\to\infty}Dx_n=Dx$.
\item[(D3)] $D$ is atomic at $0$ (see definition in Hale~\cite{hale} or
Hale and Verduyn Lunel~\cite{hale2}).
\end{itemize}
\par\smallskip From (D1) and (D2) we obtain the following
representation.
\smallskip
\begin{prop}\label{DD1D2}
 If $D\colon BU\to \R^m$ satisfies ${\rm (D1)}$ and
${\rm(D2)}$, then for each $x\in BU$
\[Dx=\int_{-\infty}^0 [d\mu(s)]\,x(s)\]
where $\mu=[\mu_{ij}]$ and $\mu_{ij}$ is a real regular Borel
measure with finite total variation
$|\mu_{ij}|(-\infty,0]<\infty$, for all $i$, $j\in\{1,\ldots,m\}$.
\end{prop}
\smallskip
\begin{proof}
>From Riesz representation theorem we obtain the above relation for
each $x$ whose components are of compact support. Moreover, if
$x\in BU$ there are an $r>0$ and a sequence of functions of
compact support $\{x_n\}_{n\in\N}\subset B_r$ with
$\n{x_n}_\infty\leq \n{x}_\infty$ such that
$x_n\stackrel{\textsf{d}\;}\to x$ as $n\to\infty$ and, from
hypothesis (D2), $\lim_{\,n\to\infty}Dx_n=Dx$. However,
\[Dx_n= \int_{-\infty}^0 [d\mu(s)]\,x_n(s)\]
and Lebesgue dominated convergence theorem yields
\[ \lim_{n\to\infty}Dx_n= \int_{-\infty}^0 [d\mu(s)]\,x(s)\,,\]
which finishes the proof.
\end{proof}
\par
Since in addition $D$ is atomic at $0$, $\det
[\mu_{ij}(\{0\})]\neq 0$, and without loss of generality, we may
assume that
\begin{equation}\label{Dformula} Dx=x(0)-\int_{-\infty}^0
[d\nu(s)]\,x(s)
\end{equation}
where $\nu=[\nu_{ij}]_{i,j\in\{1,\ldots,m\}}$, $\nu_{ij}$ is a
real regular Borel measure with finite total variation, and
$|\nu_{ij}|(\{0\})=0$ for all $i$, $j\in\{1,\ldots,m\}$. We will
denote by $|\nu|[-r,0]$ the $m\times m$ matrix
$[\,|\nu_{ij}|[-r,0]\,]$ and by $\n{\nu}_\infty[-r,0]$ the
corresponding matricial~norm.
\par
>From now on, we will assume that the operator $D$ satisfying
(D1-D3) has the form~\eqref{Dformula}. First, it is easy to check
the following result whose proof is omitted.\par\smallskip
\begin{prop}\label{existencia}
For all $h\in C([0,\infty),\R^m)$ and $\varphi\in BU$ with
$D\varphi=h(0)$, the nonhomogeneous equation
\begin{equation}\label{nohomogenea}
\left\{
\begin{array}{ll}
Dx_t=h(t)\,, & t\geq 0\,, \\
x_0=\varphi\,,
\end{array} \right.
\end{equation}
has a solution defined for all $t\geq 0$.
\end{prop}
\par
Next we obtain a bound for the solution in a finite interval
$[0,T]$, in terms of the initial data and the independent term
$h$, which in particular implies the uniqueness of the solution
of~\eqref{nohomogenea}.
\par\smallskip
\begin{lema}\label{cotaT}
Given $T>0$, there are positive constants $k_T^1$, $k_T^2$ such
that if $x$ is a solution of~\eqref{nohomogenea}, then for each
$t\in[0,T]$
\begin{equation}\label{cotaT2}
 \n{x_t}_\infty\leq k_T^1\,\sup_{0\leq
u\leq t}\n{h(u)}+ k_T^2\,\n{\varphi}_\infty\,.
\end{equation}
\end{lema}
\begin{proof}
Since $|\nu_{ij}|[-r,0]\to 0$ as $r\to 0$, for each $i$,
$j\in\{1,\ldots,m\}$, there is an $r>0$ such that
$\n{\nu}_\infty[-r,0]<1/2$. Let $x$ be a solution
of~\eqref{nohomogenea}. From~\eqref{Dformula},
\[ x(t)=h(t)+\int_{-t}^0 [d\nu(s)]\,x(t+s)+
\int_{-\infty}^{-t}[d\nu(s)]\,\varphi(t+s)\] for each $t\geq 0$.
Consequently, if $t\in[0,r]$
\[\n{x(t)}\leq \n{h(t)}+\frac{1}{2} \sup_{0\leq u\leq t}\n{x(u)}+
\n{\varphi}_\infty\n{\nu}_\infty(-\infty,0]\,,\] from which we
deduce that if $t\in[0,r]$
\begin{equation}\label{0r}
\sup_{0\leq u\leq t}\n{x(u)}\leq 2 \,\sup_{0\leq u\leq t}\n{h(u)}
+ 2\,a\,\n{\varphi}_\infty\,,
\end{equation}
where $a=\n{\nu}_\infty(-\infty,0]$. Next, let $y(t)=x(t+r)$,
which is a solution of
\[
\left\{
\begin{array}{ll}
Dy_t=h(t+r)\,, & t\geq 0\,, \\
y_0=x_r\,.
\end{array} \right.
\]
As above, we conclude that if $t\in[0,r]$
\[
\sup_{0\leq u\leq t}\n{y(u)}\leq 2 \,\sup_{0\leq u\leq
t}\n{h(u+r)} + 2\,a\,\n{x_r}_\infty\,,
\]
which together with $\n{x_r}_\infty\leq
\n{\varphi}_\infty+\sup_{0\leq u\leq r}\n{x(u)}$ and~\eqref{0r}
yields
\[
\sup_{0\leq u\leq t}\n{x(u)}\leq b \,\sup_{0\leq u\leq t}\n{h(u)}
+ c\,\n{\varphi}_\infty\,,
\]
for $t\in[r,2\,r]$ and some positive constants $b$ and $c$
independent of $h$ and $\varphi$. This way, the result is obtained
in a finite number of steps.
\end{proof}\par
Following Hale~\cite{hale}, we introduce the concept of stability
for the operator $D$. Although the initial definition is given for
the homogeneous equation, it is easy to deduce quantitative
estimates in terms of the initial data for the solution of a
non-homogeneous equation.
\begin{defi}\label{Dstable}
The linear operator $D$ is said to be {\em stable\/} if there is a
continuous function $c\in C([0,\infty),\R^+)$ with
$\lim_{\,t\to\infty}c(t)=0$ such that, for each $\varphi\in BU$
with $D\varphi=0$, the solution of the homogeneous problem
\[\left\{
\begin{array}{ll}
Dx_t=0\,, & t\geq 0 \\
x_0=\varphi\,,
\end{array} \right.
\] satisfies $\n{x(t)}\leq c(t)\,\n{\varphi}_\infty$
for each $t\geq 0$.
\end{defi}
\smallskip
\begin{prop}\label{cotacero}
 Let us assume that $D$ is stable. Then there is a
positive constant $d>0$ such that, for each $h\in
C([0,\infty),\R^m)$ with $h(0)=0$, the solution of
\[\left\{
\begin{array}{ll}
Dx_t=h(t)\,, & t\geq 0\,, \\
x_0=0\,,
\end{array} \right.
\] satisfies $\displaystyle{\n{x(t)}\leq d\sup_{0\leq u\leq t}\n{h(u)}}$
for each $t\geq 0$.
\end{prop}
\begin{proof}
Let $\{e_1,\ldots,e_m\}$ be the canonical basis of $\R^m$. A
similar proof to the one of Lemma 3.2 \S12 of Hale~\cite{hale}
shows that there are m functions $\phi_1,\ldots,\phi_m\in BU$ such
that $D\phi_j=e_j$ for each $j\in \{1,\ldots,m\}$. We will denote
by $\Phi$ the $m\times m$ matrix function
$\Phi=[\phi_1,\ldots,\phi_m]$ and by $\n{\Phi}_\infty$ the
matricial norm corresponding to the norm $\n{\cdot}_\infty$ on $BU$.
\par  Let $c\in C([0,\infty),\R^m)$ be the function given in
Definition~{\rm \ref{Dstable}}. Assume that $c$ is decreasing and
 take $T>0$ such that $c(T)<1$. From Lemma~\ref{cotaT},
$\n{x(t)}\leq k_T^1\,\sup_{\,0\leq u\leq t} \n{h(u)}$ provided
that $t\in[0,T]$. \par  If $t\geq T$, there is a $j\in \N$ such
that $t\in[j\,T,(j+1)\,T]$ and it is easy to check that
$x(t)=x^1(t-(j-1)\,T)+x^2(t-(j-1)\,T)$ where $x^1$ and $x^2$ are
the solutions of
\[
\left\{
\begin{array}{l}
Dx_t^1=0\,, \quad t\geq 0\,, \\
x_0^1=x_{(j-1)\,T}-\Phi\,h((j-1)\,T)\,,
\end{array} \right.  \quad \left\{
\begin{array}{ll}
Dx_t^2=h(t+(j-1)\,T)\,, & t\geq 0\,, \\
x_0^2=\Phi\,h((j-1)\,T)\,,
\end{array} \right.
\]
respectively. From the stability of $D$ and Lemma~\ref{cotaT}, we
deduce that
\[\n{x(t)}\leq
c(t-(j-1)\,T)\,\n{x_{(j-1)\,T}-\Phi\,h((j-1)\,T)}_\infty +
k_{2T}\,\sup_{(j-1)\,T\leq u\leq t}\n{h(u)}\,.
\]
In addition, since $t-(j-1)\,T\geq T$ and $c$ is decreasing we
conclude that
\begin{equation}\label{cotaj}
\n{x(t)}\leq c(T)\,
c_j+(c(T)\,\n{\Phi}_\infty+k_{2T})\,\sup_{0\leq u\leq
t}\n{h(u)}\,,\quad t\in[j\,T,(j+1)\,T]
\end{equation}
where $c_j=\n{x_{j\,\!T}}_\infty=\displaystyle{\sup_{0\leq u\leq
j\,\!T}\n{x(u)}}$.
\par
Let $a_T=\max\{k_T^1,\, c(T)\,\n{\Phi}_\infty+k_{2T}\}$. We have
 $c_1\leq \displaystyle{a_T\sup_{0\leq u\leq T}\n{h(u)}}$
and from~\eqref{cotaj}, if $j\geq 2$
\[c_j\leq \max\left\{c_{j-1}, c(T)\,c_{j-1}+a_T\,\sup_{0\leq u\leq
j\,T}\n{h(u)}\right\}\,.\] Hence, we check that for each $j\geq 2$
\[c_j\leq a_T\left(1+c(T)+\ldots c(T)^{j-1}\right)\sup_{0\leq u\leq
j\,T}\n{h(u)}\,,\] and again from~\eqref{cotaj} we finally deduce
that for $t\geq0$ (and hence $t\in[j\,T,(j+1)\,T]$ for some $j\geq
0$)
\[ \n{x(t)}\leq a_T\sum_{k=0}^j c(T)^k\,\sup_{0\leq u\leq t}\n{h(u)}\leq
\frac{a_T}{1-c(T)}\,\sup_{0\leq u\leq t}\n{h(u)}\,,
\]
which finishes the proof.
\end{proof}
\smallskip
\begin{teor}\label{cotanohomo}
Let us assume that $D$ is stable. Then there is a continuous
function $c\in C([0,\infty),\R^+)$ with
$\lim_{\,t\to\infty}c(t)=0$ and a positive constant $k>0$ such
that the solution of~\eqref{nohomogenea} satisfies
\[\n{x(t)}\leq c(t)\,\n{\varphi}_\infty+ k\,\sup_{0\leq
u\leq t}\n{h(u)}\]
 for each $t\geq 0$.
\end{teor}
\begin{proof}
It is not hard to check that $x(t)=x^1(t)+ x^2(t)$ where $x^1$ and
$x^2$  are the solutions of
\[
\left\{
\begin{array}{l}
Dx_t^1=\psi(t)\,h(t)\,, \quad t\geq 0\,, \\
x_0^1=\varphi\,,
\end{array} \right.  \quad \left\{
\begin{array}{ll}
Dx_t^2=(1-\psi(t))\,h(t)\,, & t\geq 0\,, \\
x_0^2=0\,,
\end{array} \right.
\]
respectively, and
\[\begin{array}{lccl} \psi\colon & [0,\infty)&
\longrightarrow &\R\\ [-.2cm] & t  &\mapsto &
\psi(t)=\begin{cases} 1-t,\; & 0\leq t\leq 1\,,\\ 0\,,\;& 1\leq
t\,.
\end{cases}
\end{array}\]
Moreover, since $y(t)=x^1(t+1)$ satisfies $ Dy_t=0$, $t\geq 0$
with $y_0=x_1^1$, the result follows from the application of
Definition~\ref{Dstable}, Proposition~\ref{cotaT} and
Proposition~\ref{cotacero} to $y$, $x^1$ on $[0,1]$, and $x^2$
respectively.
\end{proof}\par
The conclusions of Theorem~\ref{cotanohomo} are essential in what
follows. In particular, it allows us to estimate the norm of a
function $x$ in terms of the norm of the function~$(-\infty,0]\to
\R^m$, $s\mapsto Dx_s$.\par
\smallskip
\begin{prop}\label{D-1bounded}
Let us assume that $D$ is stable. Then there is a positive
constant $k>0$ such that $\n{x^h}_\infty\leq k\,\n{h}_\infty$ for
all $h\in BU$ and $x^h\in BU$ satisfying $Dx_s^h=h(s)$ for $s\leq
0$.
\end{prop}
\begin{proof} Let $x(t)$ be the solution of
\[\left\{
\begin{array}{ll}
Dx_t=h(0)\,, & t\geq 0\,, \\
x_0=x^h\,,
\end{array} \right.
\]
$\wt h(t)=\begin{cases}
h(t)\,, &t\leq 0\,,\\
h(0)\,, & t\geq 0\,,
\end{cases}$ and for $s\leq 0$ we define $y^s(t)=\begin{cases} x(t+s)\,, &
t+s\geq 0\,,\\ x^h(t+s)\,, &t+s\leq 0\,.
\end{cases}$
Then
\[\left\{
\begin{array}{ll}
Dy_t^s=\wt h(t+s)\,, & t\geq 0\,, \\
y_0^s=x^h_s\,,
\end{array} \right.
\]
and Theorem~\ref{cotanohomo} yields
\[\n{y^s(t)}\leq c(t)\,\n{x_s^h}_\infty+ k\,\sup_{0\leq u\leq t}\n{\wt h(u+s)}_\infty
\leq c(t)\,\n{x^h}_\infty+k\, \n{h}_\infty\,,\] for all $t\geq 0$
and $s\leq 0$. Hence, $\n{x^h(s)}=\n{y^{s-t}(t)}\leq
c(t)\,\n{x^h}_\infty+k\, \n{h}_\infty$,  and as $t\to\infty$ we
prove the result.
\end{proof}\par
Let $D$ be stable and given by~\eqref{Dformula}. We define the
linear operator
\begin{equation}\label{Dhat}
\begin{array}{lcclcl}
\widehat{D} \colon &BU &\longrightarrow & BU &&\\
& x & \mapsto &\widehat{D}x\colon(-\infty,0]& \to  &\R^m \\
&  & &   \hspace{1.2cm}s & \mapsto & Dx_s\,,
\end{array}
\end{equation}
that is,
$\widehat{D}x(s)=x(s)-\int_{-\infty}^0[d\nu(\theta)]\,x(\theta+s)$
for each $s\in(-\infty,0]$, which is well defined, i.e.
$h=\widehat{D}x\in BU$ provided that $x\in BU$ because $D$ is
bounded and $h(s+\tau)-h(s)=D\left(x_{s+\tau}-x_s\right)$, for all
$\tau$, $s\leq 0$. Moreover, it is easy to check that $\widehat D$
is bounded for the norm and uniformly continuous when we take the
restriction of the compact-open topology to $B_r$, i.e. given
$\varepsilon >0$ there is a $\delta(r)>0$ such that
$\textsf{d}(\widehat D x_1,\widehat{D} x_2)<\varepsilon$ for all
$x_1$, $x_2\in B_r$ with $\textsf{d}(x_1,x_2)<\delta(r)$. The next
result shows, after proving that $\widehat D$ is invertible, that
the same happens for $\widehat D^{-1}$.
\smallskip
\begin{teor}\label{inversoD}
 Let us assume that $D$ is stable. Then $\widehat{D}$ is invertible, $\widehat{D}^{-1}$
is bounded for the norm and uniformly continuous when we take the
restriction of the compact-open topology to $B_r$, i.e. given
$\varepsilon >0$ there is a $\delta(r)>0$ such that
{\upshape$\textsf{d}(\widehat{D}^{-1}h_1,\widehat{D}^{-1}h_2)<\varepsilon$}
for all $h_1$, $h_2\in B_r$ with
{\upshape$\textsf{d}(h_1,h_2)<\delta(r)$}.
\end{teor}
\begin{proof}  $\widehat{D}$ is injective
because from Proposition~\ref{D-1bounded} the only solution of
$Dx_s=0$, for $s\leq 0$ is $x=0$. To show that $\widehat{D}$ is
onto, let $h\in BU$ and $\{h_n\}_{n\in\N}\subset B_r$, for some
$r>0$, be a sequence of continuous functions whose components are
of compact support such that $h_n\stackrel{\textsf{d}\;}\to h$ as
$n\uparrow\infty$. Moreover, it is easy to choose them with the
same modulus of uniform continuity as $h$. It is not hard to check
that for each $n\in\N$ there is an $x^n\in BU$ such that
$\widehat{D}x^n=h_n$, that is, $Dx_s^n=h_n(s)$ for $s\leq 0$ and
$n\in\N$. From Proposition~\ref{D-1bounded}, $x^n\in B_{k\,\!r}$
because $\n{x^n}_\infty\leq k\,\n{h_n}_\infty\leq k\,r$ and
$\n{x^n-x^n_\tau}_\infty\leq k\,\n{h_n-(h_n)_\tau}_\infty$ for
each $\tau\leq 0$ and $n\in\N$, which implies that
$\{x_n\}_{n\in\N}$ is equicontinuous, and hence, relatively
compact for the compact open topology. Hence, there is a
convergent subsequence, let us assume the whole sequence, i.e.
there is a continuous function $x$ such that
$x^n\stackrel{\textsf{d}\;}\to x$ as $n\uparrow \infty$. From
this,  $x^n_s \stackrel{\textsf{d}\;}\to x_s$ for each $s\leq 0$
and~\eqref{Dformula} yields $Dx_s^n=h_n(s)\to Dx_s$, i.e.
$Dx_s=h(s)$ for $s\leq 0$ and $\widehat Dx=h$. It is immediate to
check that $x\in BU$ and then $\widehat{D}$ is onto, as
claimed.\par  Since $\widehat D$ is linear, bounded for the norm
and bijective, the continuity of $\widehat{D}^{-1}$ for the norm
is immediate. However, it also follows from
Proposition~\ref{D-1bounded} which reads as
$\n{\widehat{D}^{-1}h}_\infty\leq k\,\n{h}_\infty$. Finally, since
$\widehat D^{-1}$ is linear, to check the uniform continuity for
the metric on each $B_r$, it is enough to prove the continuity at
$0$, i.e. $\widehat D^{-1} x_n \stackrel{\textsf{d}\;}\to 0$ as
$n\uparrow \infty$, whenever $x_n\stackrel{\textsf{d}\;}\to 0$ as
$n\uparrow \infty$ and $\{x_n\}_{n\in\N}\subset B_r$. Let
$y^n=\widehat D^{-1} x_n\in B_{kr}$.  We extend the definition of
$y^n$ to $t\geq 0$ as the solution of
\[\left\{
\begin{array}{ll}
Dy^n_t=x_n(0)\,, & t\geq 0\,, \\
y^n_0=y^n\,.
\end{array} \right.
\]
The stability of $D$ provides
\begin{equation}\label{cotanm}
\n{y^n(t+s)}\leq c(t)\,\n{y^n}_\infty+k\,\sup_{s\leq u\leq
0}\n{x_n(u)}
\end{equation}
for all $t\geq 0$ and $s\leq 0$. Now we check that
$\{y^n\}_{n\in\N}$ converges uniformly to $0$ on each compact set
$K=[-a,0]$. Given $\varepsilon>0$, there is  a $t_0>0$ such that
$c(t_0)\,\n{y^n}_\infty< \varepsilon/2$ for each $n\in\N$.
Moreover, $x_n\to 0$ in $\wt K=[-a-t_0,0]$ and hence, there is an
$n_0$ such that for each $n\geq n_0$ we have $k\,\n{x_n}_{\wt
K}<\varepsilon/2$. Therefore, from~\eqref{cotanm} we deduce that
for all $u\in K=[-a,0]$ and $n\geq n_0$,
\[\n{y^n(u)}=\n{y^n(t_0+u-t_0)}<
\frac{\varepsilon}{2}+ \frac{\varepsilon}{2}=\varepsilon\,,\] that
is, $\n{y^n}_K<\varepsilon$, and  $\widehat D^{-1}
x_n=y^n\stackrel{\textsf{d}\;}\to 0$ as $n\uparrow \infty$, which
finishes the proof.
\end{proof}\par
A more systematic study of the properties of the linear operator
$\widehat D$ defined in~\eqref{Dhat} can be found in
Staffans~\cite{staf83b}. The next result provides a necessary and
sufficient condition for a continuous operator $D$ to be stable.
In particular, if $\widehat D$ is invertible and $\widehat D^{-1}$
is continuous for the restriction of the compact-open topology to
$B_r$, then $D$ is stable.
\smallskip
\begin{teor}\label{reciproco}
Let $D\colon BU\to\R^m$ be given by~\eqref{Dformula}  and let
$\widehat D$ be the linear operator in $BU$ defined
in~\eqref{Dhat}. The following statements are equivalent:
\begin{enumerate}
\item[{\rm (i)}] $D$ is stable.
\item[{\rm (ii)}] For each $r>0$ and each sequence
$\{x_n\}_{n\in\N}$ in $BU$  such that  $\n{\widehat D
x_n}_\infty\leq r$ and {\rm$\widehat D x_n
\stackrel{\textsf{d}\;}\to 0$} as $n\uparrow \infty$, $x_n(0)\to
0$ as $n\uparrow \infty$.
\end{enumerate}
\end{teor}
\begin{proof}
(i) $\Rightarrow$ (ii) is a consequence of
Theorem~\ref{inversoD}.\par  (ii) $\Rightarrow$ (i) For each $T>0$
we define $\mathcal{L}_T\colon \{\varphi\in BU\mid D\varphi=0\}\to
\R^m$, $\varphi\mapsto x(T)$, where $x$ is the solution of
\[\left\{
\begin{array}{ll}
Dx_t=0\,, & t\geq 0\,, \\
x_0=\varphi\,.
\end{array} \right.
\]
It is easy to check that $\mathcal{L}_T$ is well defined and
linear. In addition, from~\eqref{cotaT2} we deduce that
$\n{\mathcal{L}_T(\varphi)}\leq \n{x_T}_\infty\leq
k_T^2\,\n{\varphi}_\infty$, and hence it is bounded.\par Next we
check that $\n{\mathcal{L}_T}_\infty\to 0$ as $T\to\infty$, which
shows the stability of $D$ because $\n{x(T)}\leq
c(T)\,\n{\varphi}_\infty$ for $c(T)=\n{\mathcal{L}_T}_\infty$. Let
us assume, on the contrary, that there exist $\delta>0$,  a
sequence $T_n\uparrow\infty$ and a sequence
$\{\varphi_n\}_{n\in\N}$ with $\n{\varphi_n}_\infty\leq 1$ and
$D\varphi_n=0$ such that $\n{\mathcal{L}_{T_n}(\varphi_n)}\geq
\delta$ for each $n\in\N$. That is, $\n{x^n(T_n)}\geq \delta$
where $x^n$ is the solution of
\[\left\{
\begin{array}{ll}
Dx^n_t=0\,, & t\geq 0\,, \\
x_0^n=\varphi_n\,.
\end{array} \right.
\]
Therefore,
\[
\begin{cases}
D((x_{T_n}^n)_s)=D(x_{T_n+s}^n)=0 & \text{ if } s\in[-T_n,0]\,, \\
D((x_{T_n}^n)_s)=D((\varphi_n)_{T_n+s}) & \text{ if } s\leq
-T_n\,,
\end{cases}
\]
and taking  $r=\n{D}_\infty$, the sequence
$\{x_{T_n}^n\}_{n\in\N}\subset BU$ satisfies $\n{\widehat D
x_n}_\infty\leq r$ and
 $\widehat D x_{T_n}^n\stackrel{\textsf{d}\;}\to 0$ as $n\uparrow
\infty$. Consequently, $x_{T_n}^n(0)=x^n(T_n)\to 0$ as $n\uparrow
\infty$, which contradicts the fact that $\n{x^n(T_n)}\geq
\delta$, and finishes the proof.
\end{proof}
\smallskip
\begin{prop}
Let $D\colon BU\to\R^m$ be a stable operator given
by~\eqref{Dformula} and let $\widehat D$ be the linear operator in
$BU$ defined in~\eqref{Dhat}. Then
\[\begin{array}{lccl}
D^* \colon &BU &\longrightarrow & \R^m \\
& x & \mapsto &\widehat{D}^{-1}x(0)
\end{array}
\]
is also stable and satisfies {\rm (D1-D3)}.
\end{prop}
\begin{proof}
>From Theorem~\ref{inversoD}, we deduce that $D^*$ satisfies
(D1-D2). Hence as in Proposition~\ref{DD1D2}, there is  a real
regular Borel measure $\mu^*$ with finite total variation such
that $D^*x=\int_{-\infty}^0 [d\mu^*(s)]\,x(s)$ for each $x\in BU$.
 We can write $\mu^*= A\,\delta-\nu^*$ with
$A=[\mu^*_{ij}(\{0\})]$. We claim that $D^*$ is atomic at $0$,
i.e. $\det A\neq 0$. Assume on the contrary that $\det A=0$ and
let $v\in\R^m$ be a unitary vector with $A\,v=0$. For each
$\varepsilon>0$ we take $\varphi_\varepsilon\colon (-\infty,0]\to
\R$ with $\n{\varphi_\varepsilon}_\infty=\varphi_\varepsilon(0)=1$
and $\varphi_\varepsilon(s)=0$ for each
$s\in(-\infty,-\varepsilon]$. Let $x^{\,\varepsilon}\in BU$ be
defined by $x^{\,\varepsilon}(s)=\varphi_\varepsilon(s)\,v$. The
continuity of $\widehat D$~yields
\begin{equation}\label{D-1}
1=\n{x^{\,\varepsilon}}_\infty\leq c\,\n{\widehat
D^{-1}x^{\,\varepsilon}}_\infty\,.
\end{equation}
However, for each $s\in(-\infty,0]$
\[
\widehat{D}^{-1}x^{\,\varepsilon}(s)=D^*x_s^{\,\varepsilon}
=\varphi_\varepsilon(s)A\,v-\int_{-\infty}^0
[d\nu^*(\theta)]\,\varphi_\varepsilon(\theta+s)\,v\] and,
consequently,
$\n{\widehat{D}^{-1}x^{\,\varepsilon}}_\infty\leq\|\nu^*\|_\infty(-\varepsilon,0]$,
 which tends to $0$ as $\varepsilon\to 0$,
contradicts~\eqref{D-1} and shows that $D^*$ is atomic at $0$.
Finally, $D^*$ is stable as a consequence of
Theorem~\ref{reciproco}. Notice that $\mu^*$ is the inverse of the
measure $\mu$ for the convolution defining the operator $\widehat
D$.
\end{proof}
\section{Monotone neutral functional differential
equations}\label{monotone} Throughout this section, we will study
the monotone skew-product semiflow generated by a family of NFDEs
with infinite delay and stable $D$-operator. In particular, we
establish the 1-covering property of omega-limit sets under the
componentwise separating property and uniform stability, as in
Jiang and Zhao~\cite{jizh} for FDEs with finite delay, and Novo
{\it et al.}~\cite{NOS2007} for infinite delay. The main tool in
the proof of the result is the transformation of the initial
family of NFDEs to a family of FDEs with infinite delay in whose
study the results of Novo {\it et al.}~\cite{NOS2007} turn out to
be useful.
\par  Let
$(\Om,\sigma,\R)$ be a minimal flow over a compact metric space
$(\Om,d)$ and denote $\sigma(t,\w)=\w{\cdot}t$ for all $\w \in\Om$ and
$t\in\R$. In $\R^m$, we take the maximum norm
$\n{v}=\max_{j=1,\ldots,m}|v_j|$ and the usual partial order
relation
\[
\begin{split}
 v\le w \quad&\Longleftrightarrow\quad v_j\leq w_j\quad \text{for}\;j=1,\ldots,m\,,\\
 v<w  \quad  &\Longleftrightarrow\quad\,v\,\leq\, w\quad\text{and}
 \quad v_j<w_j\quad\text{for some}\;j\in\{1,\ldots,m\}\,.
\end{split}
\]
As in Section~\ref{stableD}, we consider the Fr\'{e}chet space
$X=C((-\infty,0],\R^m)$ endowed with the compact-open topology,
i.e.~the topology of uniform convergence over compact subsets, and
$BU\subset X$ the Banach space of bounded and uniformly continuous
functions with the supremum norm
$\n{x}_\infty=\sup_{s\in(-\infty,0]} \n{x(s)}$.
\par
Let $D\colon BU \to \R^m$ be an  autonomous and stable linear
operator satisfying hypotheses (D1-D3) and given by
relation~\eqref{Dformula}. The subset
\[BU_D^+=\{x\in
BU\,|\;Dx_s\ge 0\;\text{ for each }\;s\in (-\infty,0]\}\] is a
positive cone in $BU$,  because  it is a nonempty closed subset
$BU_D^+\subset BU$ satisfying $BU_D^+ +BU_D^+\subset BU_D^+$,
$\R^+ BU_D^+\subset BU_D^+$ and $BU_D^+\cap(-BU_D^+)=\{0\}$. As
usual, a partial order relation on $BU$ is induced, given by
\[
\begin{split}
 x\le_D y\quad&\Longleftrightarrow\quad Dx_s\le Dy_s \quad\text{for each}\;s\in (-\infty,0]\,,\\
 x<_D y    \quad&\Longleftrightarrow\quad x \leq_D y\quad\text{and}\quad x\neq
 y\,.
\end{split}
\]
\begin{nota}
{\rm Notice that, if we denote the usual partial order of $BU$
\[x\le y\quad\Longleftrightarrow\quad x(s)\le y(s)\quad\text{for each}\;s\in
(-\infty,0]\,,\] we have that $x\le_D y$ if and only if $\widehat
D\,x\leq \widehat D\,y$, where $\widehat{D}$ is defined by
relation~\eqref{Dhat}. Although in some cases they may coincide,
this new order is different from the one given by Wu and Freedman
in~\cite{WF1991}.}
\end{nota}
\par\smallskip
We consider the family of non-autonomous NFDEs with infinite delay
and stable $D$-operator \addtocounter{equation}{1}
\begin{equation}\tag*{(\arabic{section}.\arabic{equation})$_\w$}\label{Ninfdelay}
 \frac{d}{dt}Dz_t=F(\w{\cdot}t,z_t)\,, \quad t\geq 0\,,\;\w\in\Om\,,
\end{equation}
defined by a function $F\colon\Om\times BU \to\R^m$,
$(\w,x)\mapsto F(\w,x)$ satisfying the following conditions:
\begin{itemize}
\item[(F1)] $F$ is continuous on $\Om\times BU$ and locally
Lipschitz in $x$ for the norm $\n{\cdot}_\infty$. \item[(F2)] For each
$r>0$, $F(\Om\times B_r)$ is a bounded subset of $\R^m$.
\item[(F3)] For each $r>0$, $F\colon\Om\times B_r\to\R^m$ is
continuous when we take the restriction of the compact-open
topology to $B_r$, i.e.~if $\w_n\to\w$ and
$x_n\stackrel{\textsf{d}\;}\to x$ as $n\to\infty$ with $x\in B_r$,
then $\lim_{n\to\infty}F(\w_n,x_n)=F(\w,x)$. \item[(F4)] If
$x,y\in BU$ with $x\leq_D y$ and $D_j\,x=D_j\,y$ holds for some
 $j\in\{1,\ldots,m\}$, then $F_j(\w,x)\leq F_j(\w,y)$ for each $\w\in\Om$.
\end{itemize}
>From hypothesis (F1), the standard theory of NFDEs with infinite
delay (see Wang and Wu~\cite{WW1985} and Wu~\cite{jwu1991})
assures that for each $x\in BU$ and each $\w\in\Om$ the
system~\ref{Ninfdelay} locally admits a unique solution
$z(t,\w,x)$ with initial value $x$, i.e. $z(s,\w,x)=x(s)$ for each
$s\in (-\infty,0]$. Therefore, the family~\ref{Ninfdelay} induces
a local skew-product semiflow
\begin{equation}\label{Ndelaskewcoop}
 \begin{array}{cccl}
 \tau &: \R^+\times\Om\times BU& \longrightarrow & \Om\times BU\\
& (t,\w,x) & \mapsto &(\w{\cdot}t,u(t,\w,x))\,,
\end{array}
\end{equation}
where $u(t,\w,x)\in BU$ and $u(t,\w,x)(s)=z(t+s,\w,x)$ for $s\in
(-\infty,0]$.
\par
As proved in Theorem~\ref{inversoD}, the operator $\widehat{D}$
defined by relation~\eqref{Dhat} is an isomorphism of $BU$. Hence,
the change of variable $y=\widehat{D}z$ takes~\ref{Ninfdelay} to
\addtocounter{equation}{1}
\begin{equation}\tag*{(\arabic{section}.\arabic{equation})$_\w$}\label{infdelay}
 y'(t)=G(\w{\cdot}t,y_t)\,, \quad t\geq 0\,,\;\w\in\Om\,,
\end{equation}
with $G\colon\Om\times BU \to\R^m$, $(\w,x)\mapsto
G(\w,x)=F(\w,\widehat{D}^{-1}x)$ satisfying the following
conditions:
\begin{itemize}
\item[(H1)] $G$ is continuous on $\Om\times BU$ and locally
Lipschitz in $x$ for the norm $\n{\cdot}_\infty$. \item[(H2)] For each
$r>0$, $G(\Om\times B_r)$ is a bounded subset of $\R^m$.
\item[(H3)] For each $r>0$, $G\colon\Om\times B_r\to\R^m$ is
continuous when we take the restriction of the compact-open
topology to $B_r$, i.e.~if $\w_n\to\w$ and
$x_n\stackrel{\textsf{d}\;}\to x$ as $n\to\infty$ with $x\in B_r$,
then $\lim_{n\to\infty}G(\w_n,x_n)=G(\w,x)$. \item[(H4)] If
$x,y\in BU$ with $x\leq y$ and $x_j(0)=y_j(0)$ holds for some
 $j\in\{1,\ldots,m\}$, then $G_j(\w,x)\leq G_j(\w,y)$ for each $\w\in\Om$.
\end{itemize}
>From hypothesis (H1), the standard theory of infinite delay FDEs
(see Hino {\it et al.}~\cite{hino}) assures that for each $x\in
BU$ and each $\w\in\Om$ the system~\ref{infdelay} locally admits a
unique solution $y(t,\w,x)$ with initial value $x$,
i.e.~$y(s,\w,x)=x(s)$ for each $s\in (-\infty,0]$. Therefore, the
new family~\ref{infdelay} induces a local skew-product semiflow
\begin{equation}\label{delaskewcoop}
 \begin{array}{cccl}
 \widehat\tau &: \R^+\times\Om\times BU& \longrightarrow & \Om\times BU\\
& (t,\w,x) & \mapsto &(\w{\cdot}t,\widehat u(t,\w,x))\,,
\end{array}
\end{equation}
where $\widehat u(t,\w,x)\in BU$ and $\widehat
u(t,\w,x)(s)=y(t+s,\w,x)$ for $s\in (-\infty,0]$, and it is
related to the previous one~\eqref{Ndelaskewcoop} by
\begin{equation}\label{usombrero}
\widehat u(t,\w,x)= \widehat{D}\,u(t,\w,\widehat{D}^{-1}\,x)\,.
\end{equation}
\par
As a consequence, most of the results obtained in Novo {\it et
al.}~\cite{NOS2007} for the skew-product
semiflow~\eqref{delaskewcoop} can  now be translated
to~\eqref{Ndelaskewcoop}.
\par
>From hypotheses (F1) and (F2), each bounded solution
$z(t,\w_0,x_0)$ provides a relatively compact trajectory, as
deduced from Proposition 4.1 of Novo {\it et al.}~\cite{NOS2007}.
\smallskip
\begin{prop}\label{relacompact}
Let $z(t,\w_0,x_0)$ be a bounded solution of equation~{\rm
\ref{Ninfdelay}}$_{_{\!0}}$, that is, $r=\sup_{t\in\R
}\n{z(t,\w_0,x_0)}<\infty$. Then ${\rm
closure}_{X}\{u(t,\w_0,x_0)\mid t\geq 0\}$ is a compact subset of
$BU$ for the compact-open topology.
\end{prop}
\par\smallskip
>From hypotheses (F1), (F2) and (F3) and Proposition 4.2 and
Corollary 4.3 of Novo {\it et al.}~\cite{NOS2007} for the
skew-product semiflow~\eqref{delaskewcoop}, we can deduce the
continuity of the semiflow~\eqref{Ndelaskewcoop} restricted to
some compact subsets $K\subset \Om\times BU$ when the compact-open
topology is considered on $BU$.
\par\smallskip
\begin{prop}\label{continuidad}
 Let $K\subset \Om\times BU$ be a
compact set for the product metric topology and assume that there
is an $r>0$ such that $\tau_t(K)\subset \Om\times B_r$ for all
$t\geq 0$. Then the map
\[\begin{array}{cccl}
 \tau &: \R^+\times K& \longrightarrow & \Om\times BU \\
& (t,\w,x) & \mapsto &(\w{\cdot}t,u(t,\w,x))\,,
\end{array}
\]
is continuous when the product metric topology is considered.
\end{prop}
\par\smallskip
>From Proposition~\ref{relacompact}, when $z(t,\w_0,x_0)$ is
bounded we can define the omega-limit set of the trajectory of the
point $(\w_0,x_0)$ as
\[\mathcal{O}(\w_0,x_0)=\{(\w,x)\in\Om\times BU\mid \exists \,t_n\uparrow \infty\;\text{ with }
 \w_0{\cdot}t_n\to\w\,,\;
u(t_n,\w_0,x_0)\stackrel{\textsf{d}\;}\to x\}\,.\]  Notice that
the omega-limit set of a pair $(\w_0,x_0)\in\Om\times BU$ makes
sense whenever ${\rm closure}_{X}\{u(t,\w_0,x_0)\mid t\geq 0\}$ is
a compact set, because then $\{u(t,\w_0,x_0)(0)=z(t,\w_0,x_0)\mid
t\geq 0\}$ is a bounded set. Proposition~\ref{continuidad} implies
that the restriction of the semiflow~\eqref{Ndelaskewcoop} to
$\mathcal{O}(\w_0,x_0)$ is continuous for the compact-open
topology. The following result is a consequence of Proposition 4.4
 of Novo {\it et al.}~\cite{NOS2007}.
 \smallskip
\begin{prop}\label{flowexten}
Let $(\w_0,x_0)\in\Om\times BU$ be such that $\sup_{t\geq
0}\n{z(t,\w_0,x_0)}<\infty$. Then $K=\mathcal{O}(\w_0,x_0)$ is a
positively invariant compact subset admitting a flow extension.
\end{prop}
\par\noindent
>From hypothesis (F4), the monotone character of the
semiflow~\eqref{Ndelaskewcoop} is deduced.
\smallskip
\begin{prop}\label{semimono} For all $\w\in\Om$ and $x, y\in
BU$ such that $x\leq_D y$ it holds that
\[u(t,\w,x)\leq_D u(t,\w,y)\]
whenever they are defined.
\end{prop}
\begin{proof}
>From $x\leq_D y$ we know that $\widehat{D}\,x\leq \widehat{D}\,y$,
and  since (F4) $\Rightarrow$ (H4), from Proposition 4.5 of Novo
{\it et al.}~\cite{{NOS2007}} we deduce that $\widehat
u(t,\w,\widehat{D}\,x)\leq \widehat u(t,\w,\widehat{D}\,y)$
whenever they are defined, that is,
\[u(t,\w,x)=\widehat{D}^{-1}\,\widehat u(t,\w,\widehat{D}\,x)\leq_D
\widehat{D}^{-1} \widehat u(t,\w,\widehat{D}\,y)= u(t,\w,y)\,,\]
as stated.
\end{proof}\par
We establish the 1-covering property of omega-limit sets when in
addition to hypotheses (F1-F4) the componentwise separating
property and uniform stability are assumed:
\begin{itemize}
\item[(F5)]  If $x,y\in BU$ with $x\leq_D y$ and $D_i\,x<D_i\,y$
holds for some
 $i\in\{1,\ldots,m\}$, then $D_i\,z_t(\w,x)<D_i\,z_t(\w,y)$ for all $t\geq
0$ and $\w\in\Om$.
\item[(F6)] There is an $r>0$ such that all the trajectories with initial
data in $\widehat D^{-1} B_r$ are uniformly stable in $\widehat
D^{-1}B_{r'}$ for each $r'>r$,  and relatively compact for the
product metric topology.
\end{itemize}
>From relation~\eqref{usombrero} we deduce that the transformed
skew-product semiflow~\eqref{delaskewcoop} satisfies:
\begin{itemize}
\item[(H5)]  If $x,z\in BU$ with $x\leq z$ and $x_i(0)<z_i(0)$
holds for some
 $i\in\{1,\ldots,m\}$, then $y_i(t,\w,x)<y_i(t,\w,z)$ for all $t\geq
0$ and $\w\in\Om$.
\item[(H6)] There is an $r>0$ such that all the trajectories with initial
data in $B_r$ are uniformly stable in $B_{r'}$ for each $r'>r$,
and relatively compact for the product metric topology.
\end{itemize}
Finally, from Theorem 5.3 of Novo {\it et al.}~\cite{NOS2007}
applied to the skew-product semiflow~\eqref{delaskewcoop}
satisfying hypotheses (H1-H6), we obtain the next result for NFDEs
with infinite delay.
\smallskip
\begin{teor}\label{copies} Assume that hypotheses {\rm(F1-F6)} hold and let
$(\w_0,x_0)\in \Om\times \widehat D^{-1} B_r$ be such that
$K=\mathcal{O}(\w_0,x_0)\subset \Om\times \widehat D^{-1} B_r$.
Then $K=\{(\w,c(\w))\mid \w\in\Om\}$ is a copy of the base and
{\rm \[\lim_{t\to\infty}
\textsf{d}(u(t,\w_0,x_0),c(\w_0{\cdot}t))=0\,,\] {\it where}} $c:\Om\to
BU$ is a continuous equilibrium, i.e. $c(\w{\cdot}t)=u(t,\w,c(\w))$ for
any $\w\in\Om$, $t\ge 0$, and it is continuous for the
compact-open topology on $BU$.
\end{teor}
\smallskip
\begin{nota}
{\rm It is easy to check that it is enough to ask for property
(F5) (and for (H5) in the case of FDEs with infinite delay)  for
initial data in $BU$ whose trajectories are globally defined
on~$\R$.}
\end{nota}
\section{Compartmental systems}\label{compartmental}
We consider compartmental models for the mathematical description
of processes in which the transport of material between
compartments takes a non-negligible length of time, and each
compartment produces or swallows material. We provide a
non-autonomous version, without strong monotonicity assumptions,
of previous autonomous results by Wu and Freedman~\cite{WF1991}
and Wu~\cite{WW1985}.
\par
Firstly, we introduce the model with which we are going to deal as
well as some notation. Let us suppose that we have a system formed
by $m$ compartments $C_1,\ldots,C_m$, denote by $C_0$ the
environment surrounding the system, and by $z_i(t)$ the amount of
material within compartment $C_i$ at time $t$  for each
$i\in\{1,\ldots,m\}$. Material flows from compartment $C_j$ into
compartment $C_i$ through a pipe $P_{ij}$  having a transit time
distribution given by a positive regular Borel measure $\mu_{ij}$
with finite total variation  $\mu_{ij}(-\infty,0]=1$, for each
$i$, $j\in\{1,\ldots,m\}$.  Let $\wt
g_{ij}:\mathbb{R}\times\mathbb{R}\to\mathbb{R}$ be the so-called
{\em transport function\/} determining the volume of material
flowing from $C_j$ to $C_i$ given in terms of the time $t$ and the
value of $z_j(t)$ for $i\in\{0,\ldots,m\}, j\in\{1,\ldots,m\}$.
For each $i\in\{1,\ldots,m\}$, we will assume that there exists an
incoming flow of material $\tilde I_i$ from the environment into
compartment $C_i$ which only depends on time. For each
$i\in\{1,\ldots,m\}$, at time time $t\ge 0$, the compartment $C_i$
produces material itself at a rate $\sum_{j=1}^m \int_{-\infty}^0
z'_j(t+s)\,d\nu_{ij}(s)$, where $\nu_{ij}$ is a positive regular
Borel measure with finite total variation
$\nu_{ij}(-\infty,0]<\infty$ and $\nu_{ij}(\{0\})=0$, for all $i$,
$j\in\{1,\ldots,m\}$. \par Once  the destruction and creation of
material is taken into account, the change of the amount of
material of any compartment $C_i$, $1\leq i\leq m$, equals the
difference between the amount of total influx into and total
outflux out of $C_i$, and we obtain a model governed by the
following system of infinite delay NFDEs:
\begin{multline}\label{neutraleq}
\frac{d}{dt}\left[z_i(t)-\sum_{j=1}^{m}\int_{-\infty}^0
z_j(t+s)\,d\nu_{ij}(s)\right]=-\wt g_{0i}(t,z_i(t))-\sum_{j=1}^m
\wt g_{ji}(t,z_i(t))
\\ + \sum_{j=1}^m\int_{-\infty}^0\wt
g_{ij}(t+s,z_j(t+s)) \,d\mu_{ij}(s)+\tilde I_i(t),
\end{multline}
$i=1\ldots,m$.  For simplicity, we denote $\wt
g_{i0}:\R\times\R\to\R^+,\;(t,v)\mapsto\wt I_i(t)$ for
$i\in\{1,\ldots,m\}$ and let $\wt g=(\wt
g_{ij})_{i,j}:\R\times\R\to \R^{m(m+2)}$. We will assume that
\begin{enumerate}
\item[(C1)] $\wt g$ is $C^1$-{\em admissible}, i.e. $\wt g$ is $C^1$ in its second variable
and $\wt g$, $\frac{\partial}{\partial v}
 \wt g$ are uniformly continuous and bounded
on $\R\times \{v_0\}$ for all $v_0\in\R$; all its components are
monotone in the second variable, and $\wt g_{ij}(t,0)=0$ for
each~$t\in\R$;
\item[(C2)] $\wt g$ is a recurrent function, i.e. its {\em hull} is
minimal;
\item[(C3)] $\mu_{ij}(-\infty,0]=1$ and $\int_{-\infty}^0
|s|\,d\mu_{ij}(s)<\infty$;
\item[(C4)] $\nu_{ij}(\{0\})=0$ and $\sum_{j=1}^m \nu_{ij}(-\infty,0]<1$, which implies that
the operator $D\colon BU\to \R^m$, with
$D_i\,x=x_i(0)-\sum_{j=1}^m \int_{-\infty}^0
x_j(s)\,d\nu_{ij}(s)$, $i=1,\ldots,m$, is stable and satisfies
(D1-D3);
\item[(C5)]  the measures $d\eta_{ij}=c_{ij}\,d\mu_{ij}-\sum_{k=0}^m
d_{ki}\,d\nu_{ij}$ are positive, where
\[c_{ij}=\inf_{(t,\,v)\in\R^2}\frac{\partial \wt
g_{ij}}{\partial v}(t,v)\,,\; \text{ and }\;
d_{ij}=\sup_{{(t,\,v)\in\R^2}}\frac{\partial \wt g_{ij}}{\partial
v}(t,v)\,.\]
\end{enumerate}
\par In  practical cases, in which the solutions with physical interest
belong to the positive cone and the functions $g_{ij}$ are only
defined on $\R\times \R^+$, we can extend them to $\R\times\R$ by
$g_{ij}(t,-v)=-g_{ij}(t,v)$ for all $v\in\R^+$. Note that (C5) is
a condition for controlling the material produced in the
compartments in terms of the material transported through the
pipes.\par  The above formulation includes some particular
interesting cases. When the measures $\nu_{ij}$ and $\mu_{ij}$ are
concentrated on a compact set, then~\eqref{neutraleq} is a NFDE
with finite delay. When the measures $\nu_{ij}\equiv 0$,
then~\eqref{neutraleq} is a family of FDE with finite or infinite
delay.
\par
As usual, we include the non-autonomous system~\eqref{neutraleq}
into a family of non-autonomous NFDEs with infinite delay and
stable $D$-operator of the form~\ref{Ninfdelay} as follows.
\par
Let $\Om$ be the {\em hull\/} of $\wt g$, namely, the closure of
the set of mappings $\{\wt g_t\mid t\in \R\}$, with $\wt
g_t(s,v)=\wt g(t+s,v)$, $(s,v)\in \R^2$, with the topology of
uniform convergence on compact sets, which from (C1) is a compact
metric space  (more precisely from the admissibility of $\wt g$,
see Hino {\it et al.}~\cite{hino}).  Let $(\Om,\sigma,\R)$ be the
continuous flow~defined on $\Om$ by translation, $\sigma:\R\times
\Om\to \Om$, $(t,\w)\mapsto\w{\cdot}t$, with $\w{\cdot}t(s,v)= \w(t+s,v)$. By
hypothesis (C2), the flow $(\Om,\sigma,\R)$ is minimal. In
addition, if $\wt g$ is almost periodic (resp. almost automorphic)
the flow will be almost periodic (resp. almost automorphic).
Notice that these two cases are included in our formulation.
\par  Let $g:\Om\times\R\to\R^{m(m+2)},\;
(\w,v)\mapsto\w(0,v)$, continuous on $\Om\times \R$ and denote
$g=(g_{ij})_{i,j}$. It is easy to check that, for all
$\w=(\w_{ij})_{i,j}\in\Om$ and all $i\in\{1,\ldots,m\}$, $\w_{i0}$
is a function dependent only on $t$; thus, we can define
$I_i=\w_{i0},\;i\in\{1,\ldots,m\}$. Let $F:\Om\times BU\to\R^m$ be
the map defined by
\[
F_i(\w,x)=-g_{0i}(\w,x_i(0))-\sum_{j=1}^m
g_{ji}(\w,x_i(0))+\sum_{j=1}^m\int_{-\infty}^0g_{ij}(\w{\cdot}s
,x_j(s))\, d\mu_{ij}(s)+I_i(\w),
\]
for $(\w,x)\in\Om\times BU$ and $i\in\{1,\ldots,m\}$. Hence, the
family \addtocounter{equation}{1}
\begin{equation}\tag*{(\arabic{section}.\arabic{equation})$_\w$}\label{NSinfdelay}
 \frac{d}{dt}Dz_t=F(\w{\cdot}t,z_t)\,, \quad t\geq 0\,,\;\w\in\Om\,,
\end{equation}
where the stable operator $D$ is defined in (C4) and satisfies
(D1-D3), includes system~\eqref{neutraleq} when $\w=\wt g$.
\par
It is easy to check that this family satisfies hypotheses (F1-F3).
The following lemma will be useful when proving (F4) and (F5). We
omit its proof, which is analogous to the one given in Wu and
Freedman~\cite{WF1991} for the autonomous case with finite delay.
\begin{lema}\label{lemadesi}
For all $\w\in \Om$, $x,\,y\in BU$ with $x\leq_D y$ and
$i=1,\ldots,m$
\begin{equation}\label{desiorder}
F_i(\w,y)- F_i(\w,x) \geq -\sum_{j=0}^m d_{ji}\,[D_iy-D_ix]
+\sum_{j=1}^m\int_{-\infty}^0 (y_j(s)-x_j(s))\,d\eta_{ij}(s)
\end{equation}
where the measures $\eta_{ij}$ are defined in {\rm (C5)}.
\end{lema}
\par\smallskip
Condition (C5) is essential to prove the monotone character of the
semiflow. It can be improved in some cases (see Arino and
Bourad~\cite{AB1990} for the scalar one).
\smallskip
\begin{prop}\label{f4f5}
Under assumptions {\rm (C1-C5)}, the family~{\rm \ref{NSinfdelay}}
satisfies hypotheses {\rm (F4)}, {\rm (F5)} and $\Om\times BU_D^+$
is positively invariant.
\end{prop}
\begin{proof}
 Let  $x,\,y\in BU$ with $x\leq_D y$ and $D_i\,x=D_i\,y$ for some
$i\in\{1,\ldots,m\}$. From (C4), apart from the stability of the
operator $D$, it is easy to prove that the inverse operator of
$\widehat D$ defined by~\eqref{Dhat} is positive. Hence, from
$x\leq_D y$, that is, $\widehat{D}\,x\leq \widehat{D}\,y$, we also
deduce that $x\leq y$ which, together with $D_i\,x=D_i\,y$,
relation~\eqref{desiorder} and hypothesis (C5), yields
$F_i(\w,y)\geq F_i(\w,x)$, that is, hypothesis (F4) holds.
\par
Next,  we check hypothesis (F5). Let $x,\,y\in BU$ with $x\leq_D
y$ and $D_i\,x<D_i\,y$ for some
 $i\in\{1,\ldots,m\}$. Since (F4) holds, from Proposition~\ref{semimono}
 $u(t,\w,x)\leq_D u(t,\w,y)$ and as before, we deduce in this case
 that $u(t,\w,x)\leq u(t,\w,y)$, i.e $z_t(\w,x)\leq z_t(\w,y)$ for
 all $t\geq 0$ and $\w \in\Om$. Let $h(t)=D_i z_t(\w,y)- D_iz_t(\w,x)$.
 From equation~\ref{NSinfdelay} and Lemma~\ref{lemadesi}
\begin{align*}
h'(t)=& F_i(\w{\cdot}t,z_t(\w,y))-F_i(\w{\cdot}t,z_t(\w,x))\\ & \geq
-\sum_{j=0}^m d_{ji}h(t) + \sum_{j=1}^m\int_{-\infty}^0
(z_j(t+s,\w,y)-z_j(t+s,\w,x))\,d\eta_{ij}(s)\,,
\end{align*}
and again from hypothesis (C5) we deduce that $h'(t)\geq -d h(t)$
for some $d\geq 0$, which together with $h(0)>0$ yields $h(t)=D_i
z_t(\w,y)- D_iz_t(\w,x)>0$ for each $t\geq 0$ and (F5) holds.
Finally,  since $I_i(\w)\geq 0$ for each $\w\in\Om$ and
$i\in\{1,\ldots,m\}$,  and the semiflow is monotone, a comparison
argument shows that $\Om\times BU_D^+$ is positively invariant, as
stated.
\end{proof}
\par Next we will study some cases in which hypotheses
{\rm (F6)} is satisfied. In order to do this, we define $M\colon
\Om\times BU\to \R$, the {\em total mass\/}  of the system~{\rm
\ref{NSinfdelay}} as
\begin{equation}\label{mass}
M(\w,x)=\sum_{i=1}^m D_ix+ \sum_{i=1}^m \sum_{j=1}^m
\int_{-\infty}^0 \left( \int_s^0g_{ji}(\w{\cdot} \tau,x_i(\tau))\,d\tau
\right)d\mu_{ji}(s)\,,
\end{equation}
for all $\w\in\Om$ and $x\in BU$, which is well defined from
condition (C3). The  next result shows the continuity properties
of $M$ and its variation along the flow.
\begin{prop}
The total mass $M$ is a continuous function on all the sets of the
form $\Om\times B_r$ with $r>0$ for the product metric topology.
Moreover, for each $t\geq 0$
\begin{equation}\label{variacionM}
 \frac{d}{dt}M(\tau_t(\w,x))=\sum_{i=1}^m\left[I_i(\w{\cdot}
t)-g_{0i}(\w{\cdot}t,z_i(t,\w,x))\right]\,.
\end{equation}
\end{prop}
\begin{proof} The continuity follows from (D2), (C1) and (C3).
A straightforward computation similar to the one given in Wu and
Freedman~\cite{WF1991} shows that
\begin{equation}\label{variacM}
M(\w{\cdot}t,z_t(\w,x))=M(\w,x)+\sum_{i=1}^m \int_0^t
\left[I_i(\w{\cdot}s)-g_{0i}(\w{\cdot}s,z_i(s,\w,x))\right]\,ds\,,
\end{equation}
from which~\eqref{variacionM} is deduced.
\end{proof}
\par
The following lemma is essential in the proof of the stability of
solutions.
\smallskip
\begin{lema}\label{lemadis}
Let $x$, $y\in BU$ with $x\leq_Dy$. Then
\begin{equation*}
0\leq D_iz_t(\w,y)-D_i z_t(\w,x)\leq M(\w,y)-M(\w,x)
\end{equation*}
for each $i=1,\ldots,m$ and whenever $z(t,\w,x)$ and $z(t,\w,y)$
are defined.
\end{lema}
\begin{proof}
>From Propositions~\ref{f4f5} and~\ref{semimono} the skew-product
semiflow induced by~\ref{NSinfdelay} is monotone. Hence, if
$x\leq_Dy$ then $u(t,\w,x)\leq_D u(t,\w,y)$ whenever they are
defined. From this, as before, since $\widehat{D}^{-1}$ is
positive we also deduce that $x\leq y$ and $u(t,\w,x)\leq
u(t,\w,y)$. Therefore, $D_i z_t(\w,x)\leq D_i z_t(\w,y)$ and
$z_i(t,\w,x)\leq z_j(t,\w,y)$ for each $i=1,\ldots,m$. In
addition, the monotonicity of  transport functions yields
$g_{ij}(\w,z_j(t,\w,x))\leq g_{ij}(\w,z_j(t,\w,y))$ for each
$\w\in \Om$. From all these inequalities,~\eqref{mass}
and~\eqref{variacM} we deduce that
\begin{multline*}
0\leq D_iz_t(\w,y)-D_iz_t(\w,x)\leq \sum_{i=1}^m
\left[D_iz_t(\w,y) -D_iz_t(\w,x)\right]\\ \leq M(\w{\cdot}t,z_t(\w,y))-
M(\w{\cdot}t,z_t(\w,x)) \leq M(\w,y)- M(\w,x)\,,
\end{multline*}
as stated.
\end{proof}
\smallskip
\begin{prop}\label{estabilidad}
Fix $r>0$. Then given $\varepsilon>0$ there exists $\delta>0$ such
that if
 $x$, $y\in B_r$ with {\rm $\di(x,y)<\delta$} then
 $\|z(t,\w,x)-z(t,\w,y)\|\leq \varepsilon$ whenever they are
 defined.
\end{prop}
\begin{proof}
 Let $c={\displaystyle\max_i} \sum_{j=1}^m
\nu_{ij}(-\infty,0]<1$.  From the continuity of $M$, given
$\varepsilon_0=\varepsilon\,(1-c)>0$ there exists
$0<\delta<\varepsilon_0$, such that if
 $x$, $y\in B_r$ with $\di(x,y)<\delta$ then
 $|M(\w,y)-M(\w,x)|<\varepsilon_0$. Therefore, if $x,\,y\in B_r$ and $x\leq_D y$ from
 Lemma~\ref{lemadis} we deduce that $0\leq D_iz_t(\w,y)-D_i
z_t(\w,x)<\varepsilon_0$ whenever $\di(x,y)<\delta$. The
definition of $D_i$ yields \begin{align*} 0\leq
z_i(t,\w,y)-z_i(t,\w,x) & <
 \varepsilon_0+\sum_{j=1}^m \int_{-\infty}^0
[z_j(t+s,\w,y)-z_j(t+s,\w,x)]\,d\nu_{ij}(s)\\
& \leq \varepsilon_0 + \n{z_t(\w,y)-z_t(\w,x)}_\infty \sum_{j=1}^m
\nu_{ij}(-\infty,0]\,, \end{align*} from which we deduce that
$\n{z_t(\w,y)-z_t(\w,x)}_\infty(1-c)<\varepsilon_0=\varepsilon(1-c)$,
that is,  $\|z(t,\w,x)-z(t,\w,y)\|\leq \varepsilon$ whenever they
are  defined. The case in which $x$ and $y$ are not ordered
follows easily from this one.
\end{proof}\par
As a consequence, from the existence of a bounded solution for one
of the systems of the family, the boundedness of all solutions is
inferred, and this is the case in which hypothesis (F6) holds.
\smallskip
\begin{teor}\label{copiescomp} Under Assumptions {\rm (C1-C5)},
if there exists $\w_0\in\Om$ such that
{\rm\ref{NSinfdelay}}$_{_{\!0}}$ has a bounded solution, then all
solutions of {\rm\ref{NSinfdelay}} are bounded as well, hypothesis
{\rm (F6)} holds, and all omega-limit sets are copies of the base.
\end{teor}
\begin{proof}
The boundedness of all solutions is an easy consequence of the
previous proposition and the continuity of the semiflow. Let
$(\w,x)\in\Om\times BU$ and $r'>0$ such that $z_t(\w,x)\in B_{r'}$
for all $t\geq 0$. Then also from Proposition~\ref{estabilidad},
we deduce that given $\varepsilon>0$ there exists a $\delta>0$
such that
\[
\n{z(t+s,\w,x)-z(t,\w{\cdot}s,y)}=\n{z(t,\w{\cdot}s,z_s(\w,x))-z(t,\w{\cdot}
s,y)}<\varepsilon
\]
for all $t\geq 0$ whenever $y\in B_{r'}$ and {\rm
$\di(z_s(\w,x),y)<\delta$}, which shows the uniform stability of
the trajectories in $B_{r'}$ for each $r'>0$. Moreover, for each
$r>0$ there is an $r'>0$ such that $\widehat D^{-1} B_r\subset
B_{r'}$. Hence, hypothesis (F6) holds for all $r>0$ and
Theorem~\ref{copies} applies for all initial data, which finishes
the proof.
\end{proof}\par
Concerning the solutions of the original compartmental system, we
obtain the following result providing a non trivial generalization
of the autonomous case, in which  the asymptotically constancy of
the solutions was shown (see Wu and Freedman~\cite{WF1991}).
Although the theorem is stated in the almost periodic case,
similar conclusions are obtained changing almost periodicity for
periodicity, almost automorphy or recurrence, that is, all
solutions are asymptotically of the same type as the transport
functions.
\smallskip
\begin{teor}
Under Assumptions {\rm (C1-C5)} and in the almost periodic case,
if there is a bounded solution of~\eqref{neutraleq}, then there is
at least an almost periodic solution and all the solutions are
asymptotically almost periodic. For closed systems, i.e. $\wt
I_i\equiv 0$ and $\wt g_{0i}\equiv 0$ for each $i=1,\ldots,m$,
there are infinitely many almost periodic solutions and the rest
of them are asymptotically almost periodic.
\end{teor}
\begin{proof}
The first statement is an easy consequence of the previous
theorem. Let $\w_0=\wt g$.  The omega-limit of each solution
$z(t,\w_0,x_0)$ is a copy of the base
$\mathcal{O}(\w_0,x_0)=\{(\w,x(\w))\mid \w\in\Om\}$ and hence,
$z(t,\w_0,x(\w_0))=x(\w_0{\cdot}t)(0)$ is an almost periodic solution
of~\eqref{neutraleq} and \[
\lim_{t\to\infty}\n{z(t,\w_0,x_0)-z(t,\w_0,x(\w_0))}=0\,.\] The
statement for closed systems follows in addition
from~\eqref{variacM}, which implies that the mass is constant
along the trajectories. Hence, there are infinitely many minimal
subsets because from the definition of the mass and (C4), given
$c>0$ there is an $(\w_0,x_0)\in\Om\times BU^+_D$ such that
$M(\w_0,x_0)=c$ and hence $M(\w,x)=c$ for each $(\w,x)\in
\mathcal{O}(\w_0,x_0)$.
\end{proof}
\section{Long-term behavior of compartmental
systems}\label{long-term} This section deals with the long-term
behavior of the amount of material within the compartments of the
compartmental system~\eqref{neutraleq} satisfying hypotheses
(C1-C5). As in the previous section, the study of the minimal sets
for the corresponding  skew-product semiflow~\eqref{Ndelaskewcoop}
induced by the family~\ref{NSinfdelay} will be essential. In
addition to hypotheses (C1-C5) we will assume the following
hypothesis:
\begin{itemize}
\item[(C6)] Given $i\in\{0,\ldots,m\}$ and $j\in\{1,\ldots,m\}$ either
$\wt g_{ij}\equiv 0$ on $\R\times\R^+$ (and hence $g_{ij}\equiv 0$
on $\Om\times \R^+$), i.e. {\em there is not a pipe from
compartment $C_j$ to compartment $C_i$},  or for each $v>0$ there
is a $\delta_v>0$ such that $\wt g_{ij}(t,v)\geq \delta_v$ for all
$t\in\R$ (and hence $g_{ij}(\w,v)>0$ for all $\w\in\Om$ and
$v>0$). In this case we will say that the pipe $P_{ij}$
\emph{carries material (or that there is a pipe from compartment
$C_j$ to compartment $C_i$)}.
\end{itemize}\par
\noindent
  Let $I=\{1,\ldots,m\}$.
$\mathcal P(I)$ denotes, as usual, the set of all  subsets of $I$.
\smallskip
\begin{defi}
Let $\zeta:\mathcal P(I)\to\mathcal P(I),\;J\mapsto\cup_{j\in
J}\{i\in I \mid P_{ij}\mbox{ carries material}\}$. A subset $J$ of
$I$ is said to be \emph{irreducible} if $\zeta(J)\subset J$ and no
proper subset of $J$ has that property. The
system~\eqref{neutraleq} is {\em irreducible\/} if the whole set
$I$ is irreducible.
\end{defi}
\par
Note that $\zeta(I)\subset I$, so there is always some irreducible
subset of $I$. Irreducible sets detect the occurrence of
dynamically independent subsystems. Our next result gives a useful
property of irreducible sets with more than one element.
\smallskip
\begin{prop}\label{caminoJ}
If a subset $J$ of $I$ is irreducible, then, for all $i,j\in J$
with $i\neq j$, there exist $p\in\N$ and $i_1,\ldots,i_p\in J$
such that $P_{i_1i},\;P_{i_2i_1},\ldots,P_{i_pi_{p-1}}$ and
$P_{ji_p}$ carry material.
\end{prop}
\begin{proof}
Let us assume, on the contrary, that
$j\notin\cup_{n=1}^\infty\zeta^n(\{i\})=\wt J_i$. Then $\wt
J_i\subsetneq J$ and, obviously  $\zeta(\wt J_i)\subset\wt J_i$,
which contradicts the fact that $J$ is irreducible.
\end{proof}\par
Let $J_1,\ldots,J_k$ be all the irreducible subsets of $I$ and let
$J_0=I\setminus\cup_{l=1}^kJ_l$. These sets reflect the geometry
of the compartmental system in a good enough way as to describe
the long-term behavior of the solutions, as we will see below.\par
Let $K$ be any minimal subset of $\Om\times BU$ for the
skew-product semiflow induced by~\ref{NSinfdelay}. From
Theorem~\ref{copiescomp}, $K$ is of the form $K=\{(\w,x(\w))\mid
\w\in\Om\}$ where $x$ is a continuous map from $\Om$ into $BU$.
All of the subsequent results give qualitative information about
the long-term behavior of the solutions. Let us see that, provided
that we are working on a minimal set $K$, if there is no inflow
from the environment, then the total mass is constant on $K$, all
compartments out of an irreducible subset are empty and, in an
irreducible subset, either all compartments are empty or all are
never empty.  In particular, in any irreducible subset \emph{with
some outflow of material}, all compartments are empty.
\smallskip
\begin{teor}\label{closed}
Assume that  $\tilde I_i\equiv 0$ for each $i\in I$ and let
$K=\{(\w,x(\w))\mid \w\in\Om\}$ be a minimal subset of $\Om\times
BU$ with $K\subset\Om\times BU_D^+$, then
\begin{enumerate}
\item[{\rm (i)}] There exists $c\geq 0$ such that
$\left.M\right|_K\equiv c$.
\item[{\rm (ii)}] $x_i\equiv 0$ for each $i\in J_0$.
\item[{\rm(iii)}] If, for some $l\in\{1,\ldots,k\}$, there exists $j_l\in J_l$
such that $x_{j_l}\equiv 0$, then $x_i\equiv 0$ for each $i\in
J_l$. In particular, this happens if there is a $j_l\in J_l$ such
that there is outflow of material from $C_{j_l}$.
\end{enumerate}
\end{teor}
\begin{proof} We first suppose that the system is closed, i.e.
$\wt g_{0i}\equiv 0,\;\wt I_i\equiv 0$ for all $i\in I$, from
which we deduce $g_{0i}\equiv 0$ and $I_i\equiv 0$ for all $i\in
I$.
\par
(i) From~\eqref{variacM} the total mass $M$ is constant along the
trajectories and hence, $M(\w{\cdot}t, x(\w{\cdot}t))=M(\w,x(\w))$ for all
$t\geq 0$ and $\w\in\Om$, which together with the fact that $\Om$
is minimal and $M$ continuous, shows the statement.\par (ii) Let
$i\in J_0$. The set $\wt J_{i}=\cup_{n=1}^\infty\zeta^n(\{i\})$
satisfies $\zeta(\wt J_{i})\subset\wt J_{i}$ and hence, contains
an irreducible set $J_l$ for some $l\in\{1,\ldots,k\}$.
Consequently, there are $i_1,\ldots,i_p\in J_0$ and $j_l\in J_l$
such that $P_{j_li_p}$ carry material.\par  It is easy to prove
that there is an $r>0$ such that $\n{x(\w)}_\infty\leq r$ for each
$r>0$ such that $\n{x(\w)}_\infty\leq r$ for each $\w\in\Om$. We
define $M_l\colon \Om\times BU\to \R$, the {\em mass restricted to
$J_l$\/}  as
\begin{equation}\label{massjl}
M_l(\w,y)=\sum_{i\in J_l} D_iy+ \sum_{i, j\in J_l}
\int_{-\infty}^0 \left (\int_s^0g_{ji}(\w{\cdot}
\tau,y_i(\tau))\,d\tau\right)d\mu_{ji}(s)\,,
\end{equation}
which is continuous on $\Om\times B_r$.  From $x(\w)\geq_D 0$,
which also implies $x(\w)\geq 0$, and (C1), we have $0\leq
M_l(\w,x(\w))\leq M(\w,x(\w))=c$ for each $\w\in\Om$.\par  Since
$J_l$ is irreducible, for all $i\in J_l$ and $\w\in\Om$
\[ \frac{d}{dt} D_ix(\w{\cdot}t)= -\!\sum_{j\in J_l} g_{ji}(\w{\cdot}t,x_i(\w{\cdot}t)(0))+
\!\!\sum_{j\in J_l\cup J_0}\int_{-\infty}^0g_{ij}(\w{\cdot}(s+t)
,x_j(\w{\cdot}t)(s))\, d\mu_{ij}(s)
\]
because the rest of the terms vanish. Consequently,
\begin{equation}\label{derivadaml}
 \frac{d}{dt} M_l(\w{\cdot}t,x(\w{\cdot}t))=\sum_{i\in J_l}\sum_{j\in
J_0}\int_{-\infty}^0
g_{ij}(\w{\cdot}(s+t),x_j(\w{\cdot}t)(s))\,d\mu_{ij}(s)\geq 0
\end{equation}
for each $\w\in\Om$. We claim that $M_l(\w,x(\w))$ is constant for
each $\w\in\Om$. Assume, on the contrary, that there are $\w_1$,
$\w_2\in\Om$ such that $M_l(\w_1,x(\w_1))<M_l(\w_2,x(\w_2))$, and
let $t_n\uparrow\infty$ such that $\lim_{n\to\infty}
\w_2{\cdot}t_n=\w_1$. From~\eqref{derivadaml} we deduce that
$M_l(\w_2,x(\w_2))\leq M_l(\w_2{\cdot}t_n,x(\w_2{\cdot}t_n))$ for each
$n\in\N$ and taking limits as $t\to\infty$ we conclude that
$M_l(\w_2,x(\w_2))\leq M_l(\w_1,x(\w_1))$, a contradiction. Hence
$M_l(\w,x(\w))$ is constant and from~\eqref{derivadaml}
\begin{equation}\label{derivada0} \sum_{i\in J_l}\sum_{j\in J_0}\int_{-\infty}^0
g_{ij}(\w{\cdot}(s+t),x_j(\w{\cdot}t)(s))\,d\mu_{ij}(s)=0\,.
\end{equation}
Next we check that $x_{i_p}\equiv 0$. From~\eqref{derivada0} we
deduce that, for each $\w\in\Om$
\begin{equation}\label{integralip}
\int_{-\infty}^0
g_{j_li_p}(\w{\cdot}s,x_{i_p}(\w)(s))\,d\mu_{j_li_p}(s)=0\,.
\end{equation}
Assume that there is an $\w_0\in\Om$ such that
$x_{i_p}(\w_0)(0)>0$. Hence there is an $\varepsilon
>0$ with $x_{i_p}(\w_0)(s)>0$ for each $s\in(-\varepsilon,0]$, and
since $P_{j_li_p}$ carries material $g_{j_li_p}(\w_0{\cdot}
s,x_{i_p}(\w_0)(s))>0$ for $s\in(-\varepsilon,0]$. In addition,
from $\mu_{j_li_p}(-\infty,0]=1$ there is a $b\leq 0$ such that
$\mu_{j_li_p}(b-\varepsilon,b]>0$. Hence, denoting
$\w_0{\cdot}(-b)=\w_1$ we deduce~that
\[\int_{b-\varepsilon}^b
g_{j_li_p}(\w_1{\cdot}s,x_{i_p}(\w_1)(s))\,d\mu_{j_li_p}(s)>0\,,\] which
contradicts~\eqref{integralip} and shows that $x_{i_p}\equiv 0$,
as claimed. Since  $x(\w)\geq_D 0$, we have $D_{i_p} x(\w)\geq 0$
and from the definition of $D_{i_p}$ we deduce that
$D_{i_p}x(\w)=0$ for each $\w\in\Om$. Therefore, \[ 0=\frac{d}{dt}
D_{i_p}x(\w{\cdot}t)=\sum_{j=1}^m\int_{-\infty}^0g_{i_pj}(\w{\cdot}(t+s)
,x_j(\w{\cdot}t)(s))\, d\mu_{i_pj}(s)\,,
\]
from which $\displaystyle \int_{-\infty}^0
g_{i_pi_{p-1}}(\w{\cdot}s,x_{i_{p-1}}(\w)(s))\,d\mu_{i_pi_{p-1}}(s)=0$,
and as before $x_{i_{p-1}}\equiv 0$. In a finite number of steps
we check that $x_i\equiv 0$, as stated.\par  (iii) From
Proposition~\ref{caminoJ}, given $i$, $j_l\in J_l$ there exist
$p\in\N$ and $i_1,\ldots,i_p\in J_l$ such that $P_{i_1i}$,
$P_{i_2i_1},\ldots,P_{i_pi_{p-1}}$ and $P_{j_li_p}$ carry
material. If $x_{j_l}\equiv 0$, the same argument given in the
last part of (ii) shows that $x_i\equiv 0$, which finishes the
proof for closed systems.\par  Next we deal with the case when
$\tilde I_i\equiv 0$ for each $i\in I$ but the system is not
necessarily closed. We also have $I_i\equiv 0$ and
from~\eqref{variacionM} we deduce that the total mass $M$ is
decreasing along the trajectories. In particular,
\begin{equation}\label{derivadam}
 \frac{d}{dt} M(\w{\cdot}t,x(\w{\cdot}t))= -\sum_{i=1}^n g_{0i}(\w{\cdot}t,x_i(\w{\cdot}t)(0))\leq
0\,.
\end{equation}
Assume that there are $\w_1$, $\w_2\in\Om$ such that
$M(\w_1,x(\w_1))<M(\w_2,x(\w_2))$, and let $t_n\uparrow\infty$
such that $\lim_{n\to\infty} \w_1{\cdot}t_n=\w_2$. From
relation~\eqref{derivadam} we deduce that $ M(\w_1{\cdot}t_n,x(\w_1{\cdot}
t_n))\leq M(\w_1,x(\w_1))$ for each $n\in\N$ and taking limits as
$n\uparrow\infty$ we conclude that $M(\w_2,x(\w_2))\leq
M(\w_1,x(\w_1))$, a contradiction, which shows that $M$ is
constant on $K$, as stated in (i). Consequently, the derivative
in~\eqref{derivadam} vanishes and $g_{0i}(\w{\cdot}t,x_i(\w{\cdot}t)(0))=0$
for all $i\in I$, $\w\in\Om$ and $t\geq 0$. This means that
$z(t,\w,x(\w))=x(\w{\cdot}t)(0)$ is a solution of a closed system, and
(ii) and the first part of (iii) follow from the previous case.
\par
Finally, let $j_l\in J_l$ be such that there is outflow of
material from $C_{j_l}$, that is, $g_{0j_l}(\w,v)>0$ for all
$\w\in\Om$ and $v>0$. Moreover, as before,
$g_{0j_l}(\w,x_{j_l}(\w)(0))=0$ for each $\w\in\Om$, which implies
that $x_{j_l}\equiv 0$ and completes the proof.
\end{proof}\par
\smallskip
\begin{nota}\label{soluciones}
{\rm Notice that, concerning the solutions of the family of
systems~\ref{NSinfdelay} and hence the solutions of the original
system~\eqref{neutraleq} when $\w=\wt g$, we deduce that, in the
case of no inflow from the environment,
$\lim_{t\to\infty}z_i(t,\w,x_0)=0$ for all $i\in J_0$, $i\in J_l$
for compartments $J_l$ with some outflow, and each $x_0\geq_D 0$.}
\end{nota}
\smallskip
\begin{nota}
{\rm If  there is no inflow from the environment and for all
$l\in\{1,\ldots,k\}$ there is a $j_l\in J_l$ such that there is
outflow of material from $C_{j_l}$, then the only minimal set in
$\Om\times BU_D^+$ is $K=\{(\w,0)\mid \w\in\Om\}$ and all the
solutions $z(t,\w,x_0)$ with initial data $x_0\geq_D 0$ tend to
$0$ as $t\to\infty$.}
\end{nota}
\par
In a non-closed system, that is, a system which may have any
inflow and any outflow of material, if there exists a bounded
solution, i.e. all solutions are bounded as shown above, and an
irreducible set which has \emph{some inflow}, then, working on a
minimal set, all compartments of that irreducible set are nonempty
and there must be some outflow from the irreducible set.
\smallskip
\begin{teor}
Assume that  there exists a bounded solution of
{\upshape\eqref{neutraleq}} and let $K=\{(\w,x(\w))\mid
\w\in\Om\}$ be a minimal subset of $\Om\times BU_D^+$. If, for
some $l\in\{1,\ldots,k\}$, there is a $j_l\in J_l$ such that $\wt
I_{j_l}\neq 0$, i.e. there is some inflow into $C_{j_l}$, then
\begin{enumerate}
\item[{\rm (i)}] $x_i\not\equiv 0$ for each $i\in J_l$, and
\item[{\rm (ii)}] there is a $j\in J_l$ such that there is outflow
of material from $C_{j}$.
\end{enumerate}
\end{teor}
\begin{proof}
(i) Let us assume, on the contrary, that there is an $i\in J_l$
such that $x_i\equiv 0$. Then, since $x(\w)\geq_D 0$ we have that
$0\leq D_i x(\w)$ and from the definition of $D_i$ given in (C4)
we deduce that $D_i x(\w)=0$ for each $\w\in\Om$. Therefore,
\begin{equation}\label{derivadaxi}
0=\frac{d}{dt}
D_ix(\w{\cdot}t)=\sum_{j=1}^m\int_{-\infty}^0g_{ij}(\w{\cdot}(t+s)
,x_j(\w{\cdot}t)(s))\, d\mu_{ij}(s)+I_i(\w{\cdot}t)\,,
\end{equation}
for all $\w\in\Om$, $t\geq 0$ and, as in (ii) of
Theorem~\ref{closed}, we check that $x_{j_l }\equiv 0$. However,
since $\wt I_{j_l}\not\equiv 0$, there is an $\w_0\in\Om$ such
that $I_{j_l}(\w_0)>0 $, which contradicts~\eqref{derivadaxi} for
$\w=\w_0$, $i=j_l$ at $t=0$.\par  (ii) Assume on the contrary that
$g_{0j}\equiv 0$ for each $j\in J_l$. Then, if we
consider~\eqref{massjl} the restriction of the mass to $J_l$, we
check that
\begin{equation*}
 \frac{d}{dt} M_l(\w{\cdot}t,x(\w{\cdot}t))=\sum_{i\in J_l}\left[I_i(\w{\cdot}t)+\sum_{j\in
J_0}\int_{-\infty}^0
g_{ij}(\w{\cdot}(s+t),x_j(\w{\cdot}t)(s))\,d\mu_{ij}(s)\right]\geq 0
\end{equation*}
for all $\w\in\Om$ and $t\geq 0$. A similar argument to the one
given in (ii) of Theorem~\ref{closed} shows that $M_l(\w,x(\w))$
is constant for each $\w\in\Om$, which contradicts the fact that
the above derivative is strictly positive for $\w=\w_0$ at $t=0$,
and proves the statement.~\end{proof}
\par
\noindent Finally, we will change hypothesis (C6) to the following
one, slightly stronger.
\smallskip
\begin{itemize}
\item[(C6)*] Given $i\in\{0,\ldots,m\}$ and $j\in\{1,\ldots,m\}$ either
$\wt g_{ij}\equiv 0$ on $\R\times\R^+$ (and hence $g_{ij}\equiv 0$
on $\Om\times \R^+$), i.e. {\em there is not a pipe from
compartment $C_j$ to compartment $C_i$},  or for each $v\geq 0$
there is a $\delta_v>0$ such that $\frac{\partial}{\partial v} \wt
g_{ij}(t,v)\geq \delta_v$ for each $t\in\R$ (and hence
$\frac{\partial}{\partial v} g_{ij}(\w,v)>0$ for all $\w\in\Om$
and $v\geq 0$). In this case we will say that the pipe $P_{ij}$
\emph{carries material (or that there is a pipe from compartment
$C_j$ to compartment $C_i$)}.
\end{itemize}
In this case, we are able to prove that, if there exists a bounded
solution then all the minimal sets coincide both on irreducible
sets having \emph{some outflow} and out of irreducible sets.
Concerning the solutions of the initial compartmental
system~\eqref{neutraleq},
\[\lim_{t\to\infty} |z_i(t,x_0)-z_i(t,y_0)|=0\] for all $i\in J_0$,
$i\in J_l$ for compartments $J_l$ with some outflow, and all
$x_0$, $y_0\geq_D 0$.
\smallskip
\begin{teor}
Let us assume that  hypotheses {\upshape(C1-C5)} and
{\upshape(C6)*} hold and that there exists a bounded solution of
system {\upshape\eqref{neutraleq}}. Let $K_1=\{(\w,x(\w))\mid
\w\in\Om\}$ and $K_2=\{(\w,y(\w))\mid \w \in\Om\}$ be two minimal
subsets of $\Om\times BU_D^+$. Then \begin{itemize}
\item[{\rm (i)}] $x_i\equiv y_i$ for each $i\in J_0$.
\item[{\rm (ii)}] If, for some $l\in\{1,\ldots,k\}$, there is a $j_l\in J_l$
such that there is outflow of material from $C_{j_l}$ then
$x_i\equiv y_i$ for each $i\in J_l$.
\end{itemize}
\end{teor}
\begin{proof}
For each $i\in\{0,\ldots,m\}$ and each $j\in\{1,\ldots,m\}$ we
define $h_{ij}\colon \Om\to \R^+$~as
\[h_{ij}(\w)=\int_0^1\frac{\partial g_{ij}}{\partial v}(\w,s\,x_j(\w)(0)+
(1-s)\,y_j(\w)(0))\,ds\geq 0\,,\] and we consider the family of
monotone linear compartmental systems \addtocounter{equation}{1}
\begin{multline}\tag*{(\arabic{section}.\arabic{equation})$_\w$}\label{lineal}
\frac{d}{dt}D_i \hat z_t=-h_{0i}(\w{\cdot}t)\,\hat
z_i(t)-\sum_{j=1}^mh_{ji}(\w{\cdot}t)\,\hat z_i(t)\\+\sum_{j=1}^m
\int_{-\infty}^0h_{ij}(\w{\cdot}(s+t))\,\hat z_j(t+s)\,
d\mu_{ij}(s)\,,\quad \w\in\Om\,.
\end{multline}
satisfying the corresponding hypotheses (C1-C4) and (C6).
Moreover, (C5) for each of the systems~\ref{lineal}, follows from
\[\inf_{\w\in\Om}
h_{ij}(\w)\geq  \inf_{v\geq 0\,,\,\w\in\Om}\frac{\partial
g_{ij}}{\partial v}(\w,v)\,, \quad \sup_{\w\in\Om} h_{ij}(\w) \leq
\sup_{v\geq 0\,,\, \w\in\Om}\frac{\partial g_{ij}}{\partial
v}(\w,v)\,,\] and (C5) for~\eqref{neutraleq}. From the definition
of $h_{ij}$ and (C6)* we deduce that the irreducible sets for the
families~\ref{lineal} and~\ref{NSinfdelay} coincide. Consequently,
Theorem~\ref{closed} (see Remark~\ref{soluciones}) applies to this
case and we deduce that if $z_0\geq _D 0$ and $J_l$ is a
compartment with some outflow of material, then
\[
\lim_{t\to\infty} \hat z_i(t,\w,z_0)=0\,,\quad \text{ for each }
i\in J_0\cup J_l\,.
\]
The same happens for $z_0\leq_D 0$ because the systems are linear.
\par
Let $z(\w)=x(\w)-y(\w)$ for each $\w\in\Om$. It is easy to check
$\hat z(t,\w,z(\w))=z(\w{\cdot}t)(0)$ for all $\w\in\Om$ and $t\geq 0$ .
Moreover, we can find $z_1\leq_D 0$ and $z_0\geq_D 0$ such that
$z_1\leq_D z(\w)\leq_D z_0$ for each $\w\in\Om$. Hence, the
monotonicity of the induced skew-product semiflow and the
positivity of $\widehat D^{-1}$ yields
\[
\hat z(t,\w,z_1)\leq z(\w{\cdot}t)(0)\leq \hat z(t,\w,z_0) \,,\;\;
\text{ for all } \w\in\Om\,,\; t\geq 0\,,
\]
from which we deduce that $z_i\equiv 0$ for all $i\in J_0$, $i\in
J_l$ and (i) and (ii) follow.
\end{proof}\par
As a consequence, under the same assumptions of the previous
theorem, when for all $l\in\{1,\ldots,k\}$ there is outflow of
material from one of the compartments in $J_l$, there is a unique
minimal set $K=\{(\w,x(\w))\mid \w\in\Om\}$ in $\Om\times BU_D^+$
attracting all the solutions with initial data in $BU_D^+$, i.e.
\[\lim_{t\to\infty} \n{z(t,\w,x_0)-x(\w{\cdot}t)(0)}=0\,,\quad \text{
whenever }\; x_0\geq_D 0\,.\] Moreover, $x\not\equiv 0$ if and
only if there is some $j\in\{1,\ldots,m\}$ such that $\wt I_j\neq
0$, i.e. there is some inflow into one of the compartments $C_j$.
\par
For the next result, in addition to hypotheses (C1-C5) and (C6)*
we will assume the following hypothesis:
\smallskip
\begin{itemize}
\item[(C7)] If $K_1=\{(\w,x(\w))\mid
\w\in\Om\}$ and $K_2=\{(\w,y(\w))\mid \w \in\Om\}$ are two minimal
subsets of $\Om\times BU_D^+$ such that $x(\w)\leq_D y(\w)$ and
$D_ix(\w_0)=D_iy(\w_0)$ for some $\w_0\in\Om$ and
$i\in\{1,\ldots,m\}$, then $x(\w)=y(\w)$ for each $\w\in\Om$,
i.e.~$K_1=K_2$.
\end{itemize}
\par
Note that if $D_ix(\w_0)=D_iy(\w_0)$ holds for some $\w_0\in\Om$
and $i\in\{1,\ldots,m\}$, then from hypothesis (F5) we deduce that
it holds for each $\w\in\Om$.
\par   Hypothesis (C7) is
relevant when it applies to closed systems and it holds in many
cases studied in the literature.  A closed system satisfying (C7)
is irreducible. Systems with a unique compartment, studied by
Arino and Bourad~\cite{AB1990} and Krisztin and Wu~\cite{KW1996},
satisfy (C7). It follows from Theorem~\ref{closed} that
irreducible closed systems described by FDEs (see Arino and
Haourigui~\cite{ARI1987}) satisfy (C7). Closed systems given by
Wu~\cite{WF1991}, Wu and Freedman~\cite{WF1991} in the strongly
ordered case, also satisfy (C7). \par
\smallskip
\begin{defi}
Let $K_1=\{(\w,x(\w))\mid \w\in\Om\}$ and $K_2=\{(\w,y(\w))\mid \w
\in\Om\}$ be two minimal subsets. It is said that $K_1<_D K_2$ if
$x(\w)<_D y(\w)$ for each $\w\in\Om$.
\end{defi}
\par
Hypothesis (C7) allows us to classify the minimal subsets in terms
of the value of their total mass, as shown in the next result.
\smallskip
\begin{teor}
Assume that the system~\eqref{neutraleq} is closed (i.e. $\wt
I_i\equiv 0$ and $\wt g_{0i}\equiv 0$ for each
$i\in\{1,\ldots,m\}$), and hypotheses {\rm(C1-C5), (C6)*} and
{\rm(C7)} hold. Then for each $c>0$ there is a unique minimal
subset $K_c$ such that $\left.M\right|_{K_c}=c$. Moreover,
$K_c\subset \Om\times BU_D^+$ and $K_{c_1}<_D K_{c_2}$ whenever
$c_1<c_2$.
\end{teor}
\begin{proof}
Since the minimal subsets are copies of the base, and the total
mass~\eqref{mass} is constant along the trajectories and
increasing for the $D$-order because $\widehat D^{-1}$ is
positive, it is easy to check that given $c>0$ there is a minimal
subset $K_c\subset \Om\times BU_D^+$ such that
$\left.M\right|_{K_c}=c$.\par  Let $\widehat{D}$ be the
isomorphism of $BU$ defined by the relation~\eqref{Dhat}. For each
$x\in BU$ we define $x^+=\widehat{D}^{-1}\sup(0,\widehat{D}x)$.
Hence $0\leq_D x^+$, $x\leq_D x^+$ and if $y\in BU$ with $x\leq_D
y$ and $0\leq_D y$ then $x^+\leq_Dy$.
\par
Since the semiflow is monotone, from $x\leq_D x^+$ we deduce that
$u(t,\w,x)\leq_D u(t,\w,x^+)$. Since the system is closed,
$u(t,\w,0)=0$, and from $0\leq_D x^+$ we check that
$0\leq_Du(t,\w,x^+)$. Consequently $u(t,\w,x)^+\leq_Du(t,\w,x^+)$
for each $t\geq 0$.
\par
 Next we check that if $K=\{(\w,x(\w))\mid \w\in\Om\}$ is
minimal, the same happens for $K^+=\{(\w,x(\w)^+)\mid \w\in\Om\}$.
Since $x(\w{\cdot}t)=u(t,\w,x)$ for each $t\geq 0$, we deduce that
$x(\w{\cdot}t)^+=u(t,\w,x(\w))^+\leq_D u(t,\w,x(\w)^+)$, and the fact
that $\widehat D^{-1}$ is positive yields $x(\w{\cdot}t)^+\leq
u(t,\w,x(\w)^+)$ for each $t\geq 0$. In addition, since the total
mass~\eqref{mass} is constant along the trajectories and
increasing for the $D$-order, we deduce that $M(\w,x(\w)^+)=
M(\w{\cdot}t,u(t,\w,x(\w)^+))\geq
M(\w{\cdot}t,u(t,\w,x(\w))^+)=M(\w{\cdot}t,x(\w{\cdot}t)^+)$ for each $t\geq 0$.
Moreover, since $x(\w)^+$ is a continuous function in $\w$ and
$\Om$ is minimal, a similar argument to the one given in (ii) of
Theorem~\ref{closed} shows that $M(\w,x(\w)^+)$ is constant on
$\Om$ and, consequently,
$M(\w{\cdot}t,u(t,\w,x(\w)^+)=M(\w{\cdot}t,x(\w{\cdot}t)^+)$ for each $\w\in\Om$ and
$t\geq 0$. Hence, from~\eqref{mass}  we conclude that
\[ 0=\sum_{i=1}^m D_i(u(t,\w,x(\w)^+)-x(\w{\cdot}t)^+)\,,\]
that is, $D(u(t,\w,x(\w)^+))=D(x(\w{\cdot}t)^+)$ for each $\w\in\Om$ and
$t\geq 0$. Besides, it is easy to check that
$(\varphi_s)^+=(\varphi^+)_s$ whenever $\varphi\in BU$ and $s\leq
0$, from which we deduce that
$D((u(t,\w,x(\w)^+))_s)=D((x(\w{\cdot}t)^+)_s)$ for each $s\leq 0$,
$t\geq 0$ and $\w\in\Om$. That is,  $\widehat
D(u(t,\w,x(\w)^+))=\widehat D(x(\w{\cdot}t)^+)$ for each $t\geq 0$
 and $\w\in \Om$, and since $\widehat D$ is an isomorphism
$u(t,\w,x(\w)^+)=x(\w{\cdot}t)^+$ for each $t\geq 0$ and $\w\in\Om$,
which shows that $K^+$ is a minimal subset, as
$K_1=\{(\w,x(\w))\mid \w\in\Om\}$ and $K_2=\{(\w,y(\w))\mid \w
\in\Om\}$ be two minimal subsets such that
$\left.M\right|_{K_i}=c$ for $i=1,2$. We fix $\w\in\Om$. The
change of variable $\widehat z(t)= z(t)-y(\w{\cdot}t)$ takes
\ref{NSinfdelay} to
\begin{equation*} \frac{d}{dt}D\widehat z_t=G(\w{\cdot}t,\widehat z_t)\,, \quad
t\geq 0\,,\;\w\in\Om\,,
\end{equation*}
where $G(\w{\cdot}t,\widehat z_t)=F(\w{\cdot}t,\widehat
z_t+y(\w{\cdot}t))-F(\w{\cdot}t,y(\w{\cdot}t))$. It is not hard to check that this
is a new family of compartmental systems satisfying the
corresponding hypotheses (C1-C5) and (C6)* and
\[\widehat K=\{(\w,x(\w)-y(\w))\mid \w\in\Om\}\]
is one of its minimal subsets. As before
\[\widehat K^+=\{(\w,(x(\w)-y(\w))^+)\mid \w\in\Om\}\]
is also a minimal subset, and hence
\[K^+=\{(\w,y(\w)+(x(\w)-y(\w))^+)\mid \w\in\Om\}=\{(\w,z(\w))\mid \w\in\Om\}\]
is a minimal set for the initial family. For each $\w\in\Om$ we
have $z(\w)\geq_D y(\w)$. \par  Let us assume that $Dz(\w)\gg
Dy(\w)$ for each $\w\in\Om$, which implies that
$D((x(\w)-y(\w))^+)\gg 0$ for each $\w\in\Om$. Consequently,
$D((x(\w)-y(\w))^+_s)=D(((x(\w)-y(\w))_s)^+)=
D((x(\w{\cdot}s)-y(\w{\cdot}s))^+)\gg 0$ for each $s\leq 0$, and we deduce
that $\widehat D x(\w)>\widehat Dy(\w)$, i.e. $x(\w)>_D y(\w)$ for
each $\w\in\Om$, and $\left.M\right|_{K_1}>\left.M\right|_{K_2}$,
a contradiction. Hence, there are an $\w_0\in\Om$ and an
$i\in\{1,\ldots,m\}$ such that $D_iz(\w_0)=D_i y(\w_0)$, and
hypothesis (C7) provides that $z(\w)=y(\w)$ for each $\w\in\Om$.
i.e. $(x(\w)-y(\w))^+\equiv 0$ for each $\w\in\Om$, or
equivalently $x(\w)-y(\w)\leq_D 0$ for each $\w\in\Om$.  Finally,
as before, from $\left.M\right|_{K_1}=\left.M\right|_{K_2}$ we
conclude by contradiction that $x(\w)=y(\w)$ for each $\w\in\Om$,
and the minimal set $K_c$ is unique, as stated. The same argument
shows that $K_{c_1}<_D K_{c_2}$ whenever $c_1<c_2$ and finishes
the proof.
\end{proof}

\end{document}